\documentclass{article}
\usepackage{latexsym}
\newtheorem{theorem}{Theorem}[section]
\newtheorem{lemma}[theorem]{Lemma}

\def\be{\begin{equation}}
\def\ee{\end{equation}}
\def\bea{\begin{eqnarray}}
\def\eea{\end{eqnarray}}
\def\beas{\begin{eqnarray*}}
\def\eeas{\end{eqnarray*}}
 
\newcommand{\prf}{\noindent
         {\bf Proof.\ }}
\newcommand{\prfe}{\hspace*{\fill} $\Box$ 

\smallskip \noindent}

\newcommand{\supp}{\mathrm{supp}\,}
\newcommand{\diam}[1]{\mathrm{diam}({#1})}
\newcommand{\meas}[1]{\mathrm{meas}({#1})}
\newcommand{\CN}{{\cal N}} 
\def\g{\gamma}
\def\m{\mu}
\def\l{\lambda} 
\def\j{\jmath}
\def\phi{\varphi}
\def\e{\varepsilon}
\def\r{\rho}
\def\a{\alpha} 
\def\b{\beta}

\def\d{\delta}

\newcommand{\ovI}{\overline{I}}
\newcommand{\ovM}{\overline{M}}
\newcommand{\ovL}{\overline{L}}
\newcommand{\ovm}{\overline{m}}
\newcommand{\ovn}{\overline{n}}
\newcommand{\ovd}{\overline{d}}
\newcommand{\ovg}{\overline{g}}
\newcommand{\ovvM}{\overline{\overline{M}}}

\newcommand{\ovE}{\overline{E}}
\newcommand{\ovR}{\overline{R}}
\newcommand{\ovRni}{\overline{R}_{n,i}}
\newcommand{\ovWni}{\overline{W}_{n,i}}
\newcommand{\ovMni}{\overline{M}_{n,i}}
\newcommand{\ovRnk}{\overline{R}_{n,k}}
\newcommand{\ovWnk}{\overline{W}_{n,k}}
\newcommand{\ovMnk}{\overline{M}_{n,k}}

\newcommand{\ovRmi}{\overline{R}_{m,i}}
\newcommand{\ovWmi}{\overline{W}_{m,i}}
\newcommand{\ovMmi}{\overline{M}_{m,i}}

\newcommand{\ovW}{\overline{W}}
\newcommand{\ovla}{\overline{\lambda}}
\newcommand{\ovlap}{\overline{\dot\lambda}}
\newcommand{\ovlas}{\overline{\lambda}'}
\newcommand{\ovmu}{\overline{\mu}}
\newcommand{\ovp}{\overline{p}}

\newcommand{\ovj}{\overline{\j}}
\newcommand{\ovrho}{\overline{\rho}}
\newcommand{\ovmus}{\overline{\mu}'}
\newcommand{\ovvrho}{\overline{\overline{\rho}}}
\newcommand{\ovvp}{\overline{\overline{p}}}
\newcommand{\ovvj}{\overline{\overline{\j}}}
\newcommand{\ovfwbn}{\overline{FW}^2_{n}}
\newcommand{\ovfmbn}{\overline{FM}^2_{n}}
\newcommand{\ovfwan}{\overline{FW}^1_{n}}
\newcommand{\ovfman}{\overline{FM}^1_{n}}
\newcommand{\ovfwcn}{\overline{FW}^3_{n}}
\newcommand{\ovfmcn}{\overline{FM}^3_{n}}
\newcommand{\ovfwbni}{\overline{FW}^2_{n,i}}
\newcommand{\ovfmbni}{\overline{FM}^2_{n,i}}
\newcommand{\ovfwani}{\overline{FW}^1_{n,i}}
\newcommand{\ovfmani}{\overline{FM}^1_{n,i}}
\newcommand{\ovfwcni}{\overline{FW}^3_{n,i}}
\newcommand{\ovfmcni}{\overline{FM}^3_{n,i}}
\newcommand{\chid}{\chi_\delta}
\newcommand{\Ted}{T_{\e,\d}}
\newcommand{\Jedt}{J_{\e,\d,\tau}}
\newcommand{\Jed}{J_{\e,\d}}
\newcommand{\wtA}{\widetilde{A}}

\def\open#1{\setbox0=\hbox{$#1$}
\baselineskip = 0pt
\vbox{\hbox{\hspace*{0.4 \wd0}\tiny $\circ$}\hbox{$#1$}} 
\baselineskip = 10pt\!}

\def\opens#1{\setbox1=\hbox{${\scriptstyle #1}$}
\baselineskip = 0pt
\vbox{\hbox{\hspace*{0.4 \wd1} $\kern-0.35em {\scriptscriptstyle \circ}$}
\hbox{${\scriptstyle #1}$}} 
\baselineskip = 10pt\!}
 
\def\fn{\open{f}\,}
\def\dt{\partial_t}
\def\dx{\partial_x}
\def\dv{ \partial_v }
\def\dz{\partial_z}
\def\dr{\partial_r}

\newcommand{\NormRt}{\|R(t)-\ovR(t)\|}
\newcommand{\NormMt}{\|M(t)-\ovM(t)\|}
\newcommand{\NormWt}{\|W(t)-\ovW(t)\|}
\newcommand{\NormRs}{\|R(s)-\ovR(s)\|}
\newcommand{\NormWs}{\|W(s)-\ovW(s)\|}
\newcommand{\NormMs}{\|M(s)-\ovM(s)\|}
\newcommand{\NormMqqt}{\|M(t)-\ovvM(t)\|}

\def\R{{\rm I\kern-.1567em R}}
\def\N{{\rm I\kern-.1567em N}}

\begin{document} 
\title{Convergence of a particle-in-cell scheme for the spherically
       symmetric Vlasov-Einstein system}
\date{}
\author{Gerhard Rein\\
        Institut f\"ur Mathematik, Universit\"at Wien\\
        Strudlhofgasse 4, A-1090 Vienna, Austria\\ [0.2cm]
        Thomas Rodewis\\
        Nabburger Stra\ss  e 19, D-81737 Munich, Germany         
        }

\maketitle

\begin{abstract}
We consider spherically symmetric, asymptotically flat space-times
with a collisionless gas as matter model.
Many properties of the resulting Vlasov-Einstein system are
not yet accessible by purely analytical means. We present a discretized 
version of this system which is suitable for numerical implementation
and is based on the particle-in-cell technique. Convergence of the
resulting approximate solutions to the exact solution is proven and
error bounds are provided.
      
\end{abstract}

\section{Introduction}
\setcounter{equation}{0}

The properties of space-times filled with matter
are of considerable interest in general relativity.
The analytical and numerical feasibility as well as the predictions
of the models depend to a large extent on the choice of matter
model. A model for which considerable progress has been made is the
collisionless gas. We consider a large ensemble of
particles, all of which are assumed to have the same
rest mass, normalized to unity, and to move forward in time.
Therefore, their number density $f$ is a non-negative function
supported on the mass shell 
\[
PM := \left\{ g_{\alpha \beta} p^\alpha p^\beta = -1,\ p^0 >0 \right\},
\]
a sub-manifold of the tangent bundle $TM$ of the space-time manifold $M$
with metric $g_{\alpha \beta}$. We use coordinates $(t,x^a)$ with zero 
shift and corresponding canonical momenta $p^\alpha$;
Greek indices always run from 0 to 3, and Latin ones from 1 to 3.
On the mass shell $PM$ the variable $p^0$ becomes a function of the 
remaining variables $(t,x^a,p^b)$:
\[
p^0 = \sqrt{-g^{00}} \sqrt{1+g_{ab}p^a p^b} .
\] 
The Vlasov-Einstein system which governs the time evolution
of such a fully relativistic, self-gravitating collisionless gas 
now reads
\[
\partial_t f + \frac{p^a}{p^0} \partial_{x^a} f - \frac{1}{p^0}
\Gamma^a_{\b \g} p^\b p^\g  \partial_{p^a} f = 0,
\]
\[
G^{\a \b} = 8 \pi T^{\a \b},
\]
\[
T^{\a \b}
= \int p^\a p^\b f \,|g|^{1/2} \,\frac{dp^1 dp^2 dp^3}{-p_0}
\] 
where $\Gamma^\a_{\b \g}$ are the Christoffel symbols, 
$|g|$ denotes the determinant of the metric,
$G^{\a \b}$ the Einstein tensor, and $T^{\a \b}$ is
the energy-momentum tensor. All physical constants, 
in particular the speed of light,
are set to unity. We refer to \cite{A1} for a review
of results on this system. 
Note that the characteristic system 
\[
\dot x^a = p^a/p^0,\quad
\dot p^a = \Gamma^a_{\b \g} p^\b p^\g/p^0
\]
of the Vlasov equation are precisely the geodesic equations 
in the given metric, written in coordinate time $t$.
The characteristic flow leaves the mass shell $PM$ invariant.

Little is known about the behavior of solutions
to the Vlasov-Einstein system in full generality. 
In the present investigation we restrict
ourselves to the spherically symmetric and asymptotically flat case.
Using Schwarzschild coordinates the metric is taken to be
\[
ds^2 =-e^{2\m (t,r)} dt^2 +e^{2\l (t,r)}dr^2+ r^2 
        (d\theta^2 + \sin^2\theta d\phi^2),
\]
where $t \in \R, r \geq 0, \theta \in [0,\pi], \phi \in [0,2 \pi]$.
Asymptotic flatness is then expressed as
\be \label{bcinf}
\lim_{r\to\infty} \l(t,r) = \lim_{r\to\infty}\m(t,r) = 0, \quad t \in \R ,
\ee
and a regular centre is guaranteed by
\be \label{bc0}
\l(t,0)= 0, \quad t \in \R.
\ee
To formulate the system it is convenient to let $x=r\, (\sin \theta \cos \phi,
\sin \theta \sin \phi, \cos \theta)$ and, instead of using the corresponding
canonical momenta to coordinatize the mass shell over the space-time point
$(t,x)$ use frame components 
\[
v^a = p^a + (e^\l -1) \frac{x\cdot p}{r} \frac{x^a}{r};
\]
from now on $\cdot$ denotes the usual scalar product in $\R^3$
and $|v|$ the Euclidean length of $v\in \R^3$.
With the abbreviation
\[
E := \sqrt{1+|v|^2},\ \mbox{i.~e.},\ p^0=e^{-\m} E
\]
the Vlasov-Einstein system now takes the following form:
\be \label{v}
\dt f + e^{\m - \l}\frac{v}{E}\cdot \dx f -
\left( \dot \l \frac{x\cdot v}{r} + e^{\m - \l} \m'
E \right) \frac{x}{r} \cdot \dv f =0,
\ee
\bea 
e^{-2\l} (2 r \l' -1) +1 
&=&
8\pi r^2 \r , \label{f1}\\
e^{-2\l} (2 r \m' +1) -1 
&=& 
8\pi r^2 p, \label{f2} 
\eea
where
\bea
\r(t,r) = \r(t,x) 
&:=& 
\int E f(t,x,v)\,dv ,\label{r}\\
p(t,r) = p(t,x) 
&:=& 
\int \left(\frac{x\cdot v}{r}\right)^2
 f(t,x,v)\frac{dv}{E}. \label{p}
\eea
By $\dot{\phantom{\l}}$
and $\phantom{\m}'$ we denote the partial derivative with respect to
$t$ and $r$ respectively.
The phase space distribution function $f$ is assumed to be spherically
symmetric, i.\ e., $f(t,Ax,Av) = f(t,x,v)$
for every $A\in {\rm SO} (3)$ and $t\in \R,\ x,\ v\in \R^3$.
It should be noted that while the above system
is well-posed is is not the complete Vlasov-Einstein system:
Only the $00$ and $11$ components of the field equations are written.
The $01$ component is also used in the sequel; like the 
also non-trivial $22$ and $33$ components it follows
from the reduced system above, and it reads
\be
\dot \l = 
- 4 \pi r e^{\l + \m} \j, \label{f3}
\ee
where
\be
\j(t,r) = 
\j(t,x) := \int \frac{x\cdot v}{r} f(t,x,v) dv. \label{j}
\ee
Due to the symmetry the field equations have no radiative degrees
of freedom, and the solution is  determined by an initial condition
\[ 
f(0,x,v) = \fn(x,v).
\]
The initial data are taken to be spherically symmetric,
non-negative, continuously differentiable, compactly supported,
and such that
\[
\int_{|y| \leq r} \int \fn (y,v)\,dv\,dy < \frac{r}{2},\ r \geq 0.
\]
The latter condition rules out trapped surfaces at $t=0$.
We briefly discuss the main analytical results
on this system. In \cite{RR1} it was shown that each initial datum
as above launches a unique smooth solution for which
all derivatives which appear in the system exist and are continuous.
The solution can be extended in time as long as $\r$ remains bounded. 
For small initial data it is then shown that a global, geodesically
complete solution results which decays to flat Minkowski space
for $t\to \infty$. On the other hand, it is shown in
\cite{Ren1} that large data will lead to a singularity,
and in \cite{RRS1} it is shown that the first such singularity
will be at the centre. The main open problem is whether large data
still lead to solutions which are global in Schwarzschild time;
if so this would imply the cosmic censorship hypotheses.
This problem was investigated by numerical simulation in
\cite{RRS2}. Given $\fn$ as above the system was solved
numerically with initial datum
$A \fn$
where the amplitude $A>0$ was slowly increased from small to large
values. There was no indication that the solutions might blow up in finite
Schwarzschild time, and there existed a critical amplitude
$A_*$ such that for $A<A_*$ the solution dispersed while for
$A>A_*$ a black hole with mass $M(A)$ seemed to form.
The system exhibited exclusively the so-called
type I behavior: $\lim_{A \to A_*+} M(A) = M_* > 0$.
This is in sharp contrast to results
for the same sort of numerical experiment with different matter models
which exhibit type II behavior, $M_*=0$.
The findings of \cite{RRS2} are confirmed in \cite{OC}
so that considerable evidence seems to indicate that
cosmic censorship holds for the Vlasov-Einstein system.
One should note, however, that in \cite{Sh1,Sh2,Sh3,Sh4}
the formation of naked singularities was reported in a numerical
simulation of the Vlasov-Einstein system with axial symmetry.
On the other hand, certain numerical findings which were claimed
there for the Vlasov-Poisson system which is the Newtonian limit
of the Vlasov-Einstein system contradict known analytical results.

Given the significance of these questions it is desirable to
have a rigorous mathematical foundation for
the numerical scheme which is used, that is, to show that
the solutions of some appropriately discretized version
of the system converge to solutions of the Vlasov-Einstein system
if the parameters of the discretization vary appropriately,
and to give corresponding error bounds.
The aim of the present paper is to carry out this program
for the scheme used in \cite{RRS2}, a so-called particle-in-cell scheme.
This scheme is discussed in detail in the next section, the basic
set-up being as follows: The support of the initial datum is split
into small cells. In each cell a point is chosen which carries 
a weight representing the integral of $f$ over this cell. 
These particles are smeared out in space by a hat function.
From this approximation for $f$ approximations for the source terms
in the field equations and thus approximations for the fields
can be determined. With these we enter into the characteristic system
of the Vlasov equation and propagate the particles by one time step.
Then the process is repeated. This sort of numerical scheme is well
known in plasma physics and in astrophysics where it is used to
simulate the Vlasov-Maxwell or the Vlasov-Poisson system.
For convergence results for these systems we refer to
\cite{V4,V1,GS,V2,V3} and in particular to \cite{Sch1}
where the spherically symmetric Vlasov-Poisson system is considered.
General background on such schemes can be found in \cite{B,D}.

For the present system particular difficulties arise due to the fact
that the source terms in the field equations contribute to the fields 
also in a pointwise sense; as opposed to the Vlasov-Poisson and
Vlasov-Maxwell systems they are not necessarily integrated in space(time).
This lack of smoothing effect of the field equations which is 
apparent for example
in (\ref{f3}) causes considerable analytical as 
well as numerical complications,
and the analysis of a numerical scheme must be based on a careful analysis
of the analytical properties of the solutions, which
was initiated in \cite{RR1}, cf.\ also \cite{RH}. 
The paper proceeds as follows: In the next section we first
collect some additional information on the Vlasov-Einstein system
which is needed in the sequel. Then we state the discretized version
of our system and the main result. At this stage we only discretize
in phase space, thereby reducing the system to a system of ordinary 
differential
equations for the particle positions and other relevant quantities.
The proofs of our convergence results and error estimates are carried
out in Section~3 in a series of lemmas. In Section~4 we discuss the question
of how to discretize the system also in time. Since our numerical findings 
coincide with those reported in \cite{RRS1} and \cite{OC} we do not include 
them here. An open problem is to analyze the role which steady states of 
the system
play in explaining the observed type I behavior. This will be the topic
of a separate numerical investigation for which the present paper sets the
theoretical stage.

We conclude this introduction with some further references to the literature.
In \cite{RR2} and \cite{Ren2} it was shown that solutions of the 
Vlasov-Einstein
system converge to solutions of the Vlasov-Poisson system in the 
Newtonian limit. 
The latter system is much better understood than the Vlasov-Einstein system,
in particular, there exists an existence result for global smooth solutions
to general initial data, cf.~\cite{Pf,LP,Sch2}.
Another important feature of the Vlasov-Einstein system as stated above is that
it possesses a large family of steady states, cf.\ \cite{RR3,RR4,R2}.
All the results mentioned so far refer to the asymptotically flat case
which is characterized by the boundary condition (\ref{bcinf}) and from
a physics point of view represents an isolated system such as a galaxy
or globular cluster in an otherwise empty universe. 
There exists also a growing number of results on the cosmological case
of the Vlasov-Einstein system, and we refer to \cite{A1} for a discussion
of and references to these results. Finally we mention that the results of the
present paper constitute the major part of the second author's PhD thesis
\cite{Ro}.

\section{The semi-discretized approximation---main result}
\setcounter{equation}{0}

An initial datum $\fn$ as specified above launches a continuously
differentiable and spherically symmetric solution $f$. Let
$[0,T]$ denote any time interval on which this solution exists.
For numerical investigations it is important to make use of the
spherical symmetry in such a way as to reduce the number of independent
variables. We introduce
\[
r := |x|,\ w := \frac{x \cdot v}{r},\ L = x^2 v^2 - (x \cdot v)^2
= |x \times v|^2
\]
so that
\be \label{edef}
E=\sqrt{1 + w^2 + L/r^2}.
\ee
It can be shown that $f$ must be of the form
\[
f(t,x,v)=f(t,r,w,L).
\]
Instead of writing the Vlasov equation in these variables 
we write down its characteristic system which is equivalent
but more relevant for what follows:
\bea
\dot r 
&=&
e^{\m - \l} \frac{w}{E}, \label{cr} \\
\dot w
&=&
- \dot \l w - e^{\m - \l} \m' E 
+ e^{\m - \l} \frac{L}{r^3 E} \label{cw},\\
\dot L 
&=&
0 \label{cL}.
\eea
Note that the quantity $L$, the modulus of angular momentum squared,
is conserved along characteristics so that the
characteristic system is essentially two dimensional.
The field equations (\ref{f1}), (\ref{f2}), (\ref{f3}) remain unaffected
by the above change of variables,
and the source terms can be written as
\bea
\r(t,r) 
&=& 
\frac{\pi}{r^2} \int_{-\infty}^\infty \int_0^\infty E f(t,r,w,)
\,dL\, dw ,\label{rrwf}\\
p(t,r) 
&=& 
\frac{\pi}{r^2} \int_{-\infty}^\infty \int_0^\infty  
\frac{w^2}{E}
 f(t,r,w,L)\, dL\,dw, \label{prwf}\\ 
\j(t,r) 
&=& 
\frac{\pi}{r^2} \int_{-\infty}^\infty \int_0^\infty w f(t,r,w,L) 
\,dL\,dw. \label{jrwf}
\eea
In the sequel we will use both $(x,v)$ and $(r,w,L)$ in our
discussion, but the numerical scheme will be formulated in the
latter coordinates.
Let us denote by $(X,V)(s,t,x,v)$ or $(R,W,L)(s,t,r,w,L)$ the solution
of the characteristic system which at time $s=t$ takes the value
$(x,v) \in \R^6$ or $(r,w,L) \in [0,\infty[\times \R \times [0,\infty[$.
Then (by abuse of notation)
\[
f(t,x,v)=\fn ((X,V)(0,t,x,v)) = \fn ((R,W,L)(0,t,r,w,L))= f(t,r,w,L).
\]
It is important that the field equations can be solved explicitly
for $\l$ and $\m$ and their relevant derivatives.
Using (\ref{bc0}) the field equation (\ref{f1})
can be integrated to yield
\be \label{leq}
e^{-2\l(t,r)} = 1 - \frac{2 m(t,r)}{r}
\ee
where
\be \label{mdef}
m(t,r) := 4 \pi \int_0^r s^2 \r(t,s) \, ds;
\ee
the right hand side of (\ref{leq}) is positive initially by assumption on $\fn$
and remains so on the existence interval of the solution. 
It is worthwhile to note that $m(t,\infty)$ is a conserved quantity
of the system, the ADM mass. Another conserved quantity, related to the 
conservation of the number of particles, is
\be \label{partn}
\int\!\!\!\! \int e^{\l(t,r)}  f(t,x,v)\, dv\,dx.
\ee
Next we solve (\ref{f2}) for $\m'$,
\be \label{mupeq}
\m'(t,r) = e^{2\l(t,r)}\left( \frac{m(t,r)}{r^2} + 4 \pi r p(t,r) \right),
\ee
and using (\ref{bcinf}) this is integrated to give
\be \label{mueq}
\m(t,r) = - \int_r^\infty \m'(t,s) \,ds.
\ee
From (\ref{f1}) we also have
\be \label{lapeq}
\l'(t,r) = e^{2\l(t,r)} \left( - \frac{m(t,r)}{r^2} +  4 \pi r \r(t,r) \right)
\ee
which will become relevant shortly. Obviously, $\l \geq 0$,
$\m \leq 0$, and by adding (\ref{mupeq}) and (\ref{lapeq})
and observing the boundary condition (\ref{bcinf}), $\l + \m \leq 0$.
As noted above, $f$ is constant along
characteristics, but the characteristic flow is {\em not} volume
preserving as can be seen from the factor $e^\l$ in (\ref{partn}).
This complication is the price we pay for writing the system
in non-canonical momentum variables which make the explicit
representations of the metric quantities in terms of the source
terms possible and greatly simplify the Vlasov equation.
We need to know
how integrals of $f$ over pieces of phase space evolve:  
\begin{lemma} \label{phspvol}
Let $Z(\cdot,t,z) = (X,V)(\cdot,t,x,v)$ denote the solution of the
characteristic system of (\ref{v}) with $Z(t,t,z)=z\in \R^6$. Then
\[ 
\det \dz Z(s,t,z) = e^{\l(t,x) - \l(s,X(s,t,z))}. 
\]
Let $A \subset \R^6$ be measurable 
and $A(t) := Z(t,0,A)$. Then we have for the 
solution $f$ and any continuous function $g$,
\[ 
\int_{A(t)} g(t,z) f(t,z) \, dz
= 
\int_A g(t,Z(t,0,z)) \fn(z) e^{\l(0,r)-\l(t,R(t,0,z))}\, dz
\]
and
\[
\frac{d}{dt} \int_{A(t)} f(t,z) \, dz 
= 
- \int_{A(t)} f(t,z) \left(\dot \l(t,r) + \dot R(t,0,Z(0,t,z)) \l'(t,r)
\right) dz.
\]
\end{lemma}
\prf
The first assertion follows from Liouville's Theorem, and the rest follows by
a change of variables and straight forward computation. \prfe

Since we will not be too careful in distinguishing between a symmetric
subset $A\subset \R^3\times \R^3$ and 
the subset $\wtA$ of $(r,w,L)$-space
which describes $A$, that is
\[
A= \left\{ (x,v)\in \R^6 \mid \left(|x|, x \cdot v/|x|, |x \times v|^2\right) 
\in \wtA \right\},
\]
we note that
\[
\meas{A} =  4\pi^2\meas{\wtA}.
\]
This follows from $dv = \pi\, r^{-2} dw\,dL$ and $dx = 4 \pi r^2 dr$.

We make the following additional assumption on the initial datum $\fn$:
With respect to the $(r,w,L)$-variables,
\be \label{suppbound}
\supp \fn \subset [r_\mathrm{min},r_\mathrm{max}] \times
[w_\mathrm{min},w_\mathrm{max}] \times
[L_\mathrm{min},L_\mathrm{max}]
\ee
with constants $0<r_\mathrm{min}<r_\mathrm{max}$, 
$w_\mathrm{min}<w_\mathrm{max}$, $0<L_\mathrm{min}<L_\mathrm{max}$.
Since $L$ is conserved
along particle trajectories, the particles remain away from zero in space
as long as their momenta $v$ remain bounded.
This assumption, which was also made in \cite{Sch1}, avoids problems
arising from particles at the origin where our coordinates $(r,w,L)$
are not suitable. In \cite{Ro} it is shown for the Vlasov-Poisson 
system how this restriction can be avoided
by switching to Cartesian coordinates in the neighborhood of the origin,
but we prefer to avoid at least this technical
complication. In passing we note that the numerical
investigations in \cite{RRS2,OC} considered data satisfying this restriction.

We are now ready to formulate the \\
{\bf Semi-discretized approximation:}
We decompose the support of $\fn$ in $(r,w,L)$-space
or a supset of the form (\ref{suppbound})
into disjoint, connected, and 
measurable sets $A_n,\ n \in \{1,2,\dots,N\}$, with 
\[
\diam{A_n} \leq \e,\ \meas{A_n} \geq \frac{1}{100} \e^3,\
n \in \{1,2,\dots,N\},
\]
where the fineness $0<\e \ll 1$ is a small parameter;
the factor $1/100$ can be replaced by any positive constant 
which is then kept fixed.
In each cell $A_n$ we fix a point $(r_n,w_n,L_n) \in A_n$ 
and define 
\[
(R_n,W_n,L_n)(t):= (R,W,L)(t,0,(r_n,w_n,L_n)),
\] 
the corresponding solutions of the characteristic system 
(\ref{cr}), (\ref{cw}), (\ref{cL}).
The cell in $(x,v)$-space which corresponds to $A_n$ is denoted by the same
symbol.
Moreover, we define 
\be \label{antmntdef}
A_n(t) := Z(t,0,A_n),\quad M_n(t) := \int_{A_n(t)} f(t,z) \, dz;
\ee
note that Lemma~\ref{phspvol} tells us how
$M_n$ evolves. 

Our aim is to obtain approximations to $R_n(t),W_n(t),M_n(t)$
as solutions of an appropriately set-up system of ordinary equations.
Throughout the paper approximating quantities will be denoted
with a bar.
Suppose we have $N$ points in $(r,w,L)$-space with coordinates
\beas
\ovR &=& (\ovR_1,\ldots,\ovR_N) \in \R_+^N,\\
\ovW &=& (\ovW_1,\ldots,\ovW_N) \in \R^N,\\
\ovL &=& (\ovL_1,\ldots,\ovL_N) \in \R_+^N,
\eeas
carrying weights
\[
\ovM = (\ovM_1,\ldots,\ovM_N) \in \R_+^N.
\]
From this information we
generate approximations for the source terms at any $r\geq 0$:
Following (\ref{edef}) we abbreviate
\be \label{discredef}
\ovE_n := \ovE(\ovR_n,\ovW_n,\ovL_n):= \sqrt{1 + \ovW_n^2 + \ovL_n/\ovR_n^2}
\ee
and define
\beas
\ovrho(r,\ovR,\ovW,\ovL,\ovM) 
&:=& 
\frac{1}{4 \pi r^2 \d} \sum_{n=1}^N \ovE_n \ovM_n \chid(r-\ovR_n),\\
\ovp(r,\ovR,\ovW,\ovL,\ovM)
&:=& 
\frac{1}{4 \pi r^2 \d} \sum_{n=1}^N \frac{\ovW_n^2}{\ovE_n} 
\ovM_n \chid(r-\ovR_n),\\
\ovj(r,\ovR,\ovW,\ovL,\ovM)
&:=& 
\frac{1}{4 \pi r^2 \d} \sum_{n=1}^N \ovW_n \ovM_n \chid(r-\ovR_n).
\eeas
Here the hat function $\chid$ which is used to smear out the particles
is defined by
\[
\chid (\zeta) := 
\left\{
\begin{array}{ccl}
1 - |\zeta|/\d &,& |\zeta| \leq \d,\\
0&,& |\zeta| > \d,
\end{array}
\right.\quad \zeta \in \R,
\]
where $0<\d \ll 1$ is another small parameter;
the necessary relations between $\d$ and $\e$ will be specified shortly.
Suppressing the variables $\ovR,\ovW,\ovL,\ovM$ for the moment
we can now define approximations for the various quantities
appearing on the right hand side of the characteristic system etc.:
\beas
\ovm(r) 
&:=& 
4 \pi \int_0^r s^2 \ovrho(s) \, ds,\\
e^{-2 \ovla(r)} 
&:=& 
1 - \frac{2 \ovm(r)}{r},\\
\ovmus(r) 
&:=& 
e^{2 \ovla(r)}\left( \frac{\ovm(r)}
{r^2} + 4 \pi r \ovp(r) \right),\\
\ovmu(r) 
&:=& 
- \int_r^\infty \ovmus(s) \, ds,\\
\ovlas(r) 
&:=&
e^{2 \ovla(r)}\left(- \frac{\ovm(r)}{r^2} +  4 \pi r \ovrho(r) \right),\\
\ovlap(r) 
&:=&
- 4 \pi r e^{\ovla(r)+\ovmu(r)} \ovj(r).
\eeas
The approximations $\ovR_n(t)$, $\ovW_n(t)$, $\ovM_n(t)$ to the true
quantities $R_n(t)$, $W_n(t)$, $M_n(t)$ are now defined as the solutions of
of the following autonomous system of $3\,N$ ordinary differential equations:
\bea
\dot{\ovR}_n
&=&
e^{(\ovmu-\ovla)(\ovR_n)} \frac{\ovW_n}{\ovE_n}, \label{rapde}\\ 
\dot{\ovW}_n
&=&
e^{(\ovmu-\ovla)(\ovR_n)}
\frac{L_n}{\ovR_n^3\ovE_n}
-\ovlap(\ovR_n) \ovW_n 
- e^{(\ovmu-\ovla)(\ovR_n)} \ovmus(\ovR_n) \ovE_n,\ \label{wapde}\\
\dot{\ovM}_n
&=& 
- \left( \ovlap(\ovR_n) + 
e^{(\ovmu-\ovla)(\ovR_n)} \frac{\ovW_n}{\ovE_n} 
\ovlas(\ovR_n) \right) \ovM_n \label{mapde}
\eea
with initial conditions $(\ovR_n,\ovW_n)(0)=(r_n,w_n)$ and $\ovM_n(0)=M_n(0)$,
$n=1,\ldots,N$.
Clearly, we set $\ovL_n=L_n$.
Moreover, the dependence of $\ovla,\ \ovmu$, and their various derivatives
on the coordinates and weights
of all the other particles has been suppressed in the
notation above.

Given a solution $\ovR(t), \ovW(t), \ovM(t)$
of the above discretized system it is suggestive to abuse our notation
as follows:
\[
\ovrho(t,r):=\ovrho(r,\ovR(t),\ovW(t),\ovM(t),\ovL),
\]
and analogously for all the other source terms, metric coefficients,
and their various derivatives.
Note that $\ovlap$ is just a name. In general
$\ovlap \neq \dt \ovla$, whereas $\ovmus = \dr \ovmu$
and $\ovlas = \dr \ovla$.
Nevertheless, we will use the notation 
$\ovlap$ because we will compare $\ovlap$ to $\dot \l$ 
during the convergence proof.

Before we state the main result we specify certain 
quantities which we use to control the 
errors. By $\|\cdot\|_\infty$ we denote the $L^\infty$ norm with respect to
$r \in [0,\infty[$. For the characteristics we define
\beas
\|R(t)-\ovR(t)\| 
&:=& 
\max\bigl\{|R_n(s)-\ovR_n(s)| \mid n \in \{1,\ldots,N\},\ s \in [0,t]\bigr\},\\
\|W(t)-\ovW(t)\| 
&:=& 
\max\bigl\{|W_n(s)-\ovW_n(s)| \mid n \in \{1,\ldots,N\},\ s \in [0,t]\bigr\},\\
\|M(t)-\ovM(t)\| 
&:=& 
\e^{-3} 
\max\bigl\{|M_n(s)-\ovM_n(s)| \mid n \in \{1,\ldots,N\},\ s \in [0,t]\bigr\};
\eeas
since the weights $M_n$ are of order $\mathrm{O}(\e^3)$ we have inserted a 
suitable factor.
We are now able to state the main result:
\begin{theorem} \label{main}
Assume in addition that
$\fn \in C^2(\R^6)$ so that the solution has the same regularity
with respect to $x$ and $v$. Then there exist constants $C_1>0$ 
and $C_2>0$ depending only on the approximated solution $f$
restricted to the time interval $[0,T]$ such that the following holds:
If $0<\e\leq \d \leq C_1$ then any solution of the discretized system stated
above exists on $[0,T]$ and satisfies the following estimates for $t\in [0,T]$:
\[
\|R(t)-\ovR(t)\| + \|W(t)-\ovW(t)\| + \|M(t)-\ovM(t)\|
\leq C \left( \d^2 + \frac{\e}{\d}\right),
\]
\[
\|\lambda(t)-\ovla(t)\|_\infty + \|\mu(t)-\ovmu(t)\|_\infty
+ \|m(t)-\ovm(t)\|_\infty
\leq C \left( \d^2 + \frac{\e}{\d}\right),
\]
\[
\|\r(t)-\ovrho(t)\|_\infty + \|p(t)-\ovp(t)\|_\infty
+ \|\j(t)-\ovj(t)\|_\infty
\leq C \left( \d + \frac{\e}{\d^2}\right).
\]
\end{theorem}
Note that the error estimates for the source terms are one order worse than
the other estimates; this is no coincidence.

{\bf Remark.}
Strictly speaking there exists in the literature
no proof for the assertion that
the regularity assumption $\fn \in C^2(\R^6)$ will propagate to the 
solution. However,
$\partial_{x_i} f$ will satisfy a system of differential
equations which is obtained from the Vlasov-Einstein system by
differentiating the various equations accordingly, and this system
will be linear in $\partial_{x_i} f$ so that this derivative will
be $C^1$ on the existence interval
of $f$. Turning this into a proof would be a lengthy exercise. 
A version of our theorem without the additional regularity assumption
is also possible. However, we then no longer obtain error bounds
for the source terms, and instead of the regularity assumption
an a-priori bound of the approximating source terms, uniform in
$\e$ and $\d$  has to be made which can be monitored during a
run-time situation. Since the statement of this second, weaker form 
of our result requires more technical preparation we postpone it
to the next section.

\section{The semi-discretized approximation---proofs}
\setcounter{equation}{0}

First we need to introduce some further
information on the solutions of the Vlasov-Einstein system.
From \cite{RR1} we know that
any initial datum $\fn\in C^1(\R^6)$ as specified above launches a unique
solution $f$. Let $[0,T]$ be any  time interval on which
this solution exists. Then
$f \in C^1([0,T] \times \R^6)$ with respect to $(x,v)$ and
\[ 
\r, p, \j\in C^1([0,T] \times [0,\infty[),\ 
m, \l,\m \in C^2([0,T] \times [0,\infty[)
\]
as functions of $t$ and $r$. This solution together with the time interval
are now kept fixed, and we will need the following bounds:

\begin{lemma} \label{bounds}
There is a constant $D \geq 1$ such that for all $t \in [0,T]$ and 
$ r \geq 0$, 
\beas
\r(t,r),\ p(t,r),\ |\j(t,r)|,\ 
e^{2 \l(t,r)},\ |\dot\l(t,r)|,\
|\l'(t,r)|,&& \\
|\m(t,r)|,\ |\m'(t,r)|,\ |\dot\l'(t,r)|,\ |\l''(t,r)|,
\ |\m''(t,r)| &\leq& D, 
\eeas
and for all $t \in [0,T]$ and $(r,w,L) \in \supp \fn$,
\[
R(t,0,r,w,L),\ \frac{1}{R(t,0,r,w,L)},\ |\dot R(t,0,r,w,L)|,\ 
|W(t,0,r,w,L)| \leq D.
\]
Moreover, $D$ can be chosen such that
\[
\e^{-3} M_n(t) \leq D,\ n = 1,\ldots, N,
\]
for any discretization of $\supp \fn$ as specified above.
\end{lemma}

\prf
The right hand side of (\ref{cr}) is bounded by $1$, and hence
the source terms vanish outside the light cone, that is for all $(t,r)$ with
$r \geq r_{\max} + t$. Together with the regularity of the solution
this implies the first set of estimates. Since
$V(t,0,z)$ is uniformly bounded for  $t \in [0,T]$ and $z \in \supp{\fn}$
and $L$ is conserved we get a lower bound on $R$. Finally,
by Lemma~\ref{phspvol},
\[ 
M_n(t)=\int_{A_n(t)} f(t,z) \, dz
\leq
\int_{A_n} \fn(z) e^{\l(0,r)}\, dz \leq \|\fn e^{\l(0)}\|_\infty
\meas{A_n}, 
\] 
and the proof is complete.
\prfe 

All the functions on the right hand side of  (\ref{rapde}),
(\ref{wapde}), (\ref{mapde}) are Lipschitz-continuous in the unknowns. 
Let $[0,\Ted]$ denote the maximal interval on which the approximate
solution exists and satisfies the bounds 
\be \label{teda}
\e^{-3} \ovM_n(t), \ovR_n(t), \frac{1}{\ovR_n(t)}, |\ovW_n(t)|,
\|e^{2 \ovla(t)}\|_\infty \leq 2 D,\ n=1,\ldots,N .
\ee
Here $D$ is the constant used in Lemma \ref{bounds}. 
It must be emphasized that $\Ted$ depends on the 
particular discretization of the support of $\fn$ and not just on
$\e$ and $\d$.

We can now formulate the weaker version of our main result which does not
require the additional regularity of the approximated solution:
\begin{theorem} \label{mainw}
Let a family of approximating solutions
parameterized by $\e$ and $\d$ (or a sequence of such) be given
and assume that for some constant $C^\star>0$,
\[
\ovrho(t) \leq C^\star,\ 0<\e\leq\d\leq\frac{1}{4 D},\
t \in [0,\Ted].
\] 
Then there exist constants
$C_1>0$ and $C_2>0$ depending on the restriction of the approximated solution
$f$ to the time interval $[0,T]$ and on $C^\ast$ such that the following holds:
If $0<\e\leq \d \leq C_1$ then any solution of the discretized system stated
above exists on $[0,T]$ and satisfies the following estimates for $t\in [0,T]$:
\beas
\|R(t)-\ovR(t)\| + \|W(t)-\ovW(t)\| + \|M(t)-\ovM(t)\|
&\leq& 
C \left( \d + \frac{\e}{\d}\right),\\
\|\lambda(t)-\ovla(t)\|_\infty + \|\mu(t)-\ovmu(t)\|_\infty
+ \|m(t)-\ovm(t)\|_\infty
&\leq C&
\left( \d + \frac{\e}{\d}\right).
\eeas
\end{theorem}
Note that in a run-time situation the assumption on the
boundedness of $\ovrho$ can be monitored
by looking at the numerical data. 

Throughout the rest of this section we fix discretization parameters
\[
0<\e\leq \d \leq \frac{1}{4D}.
\]
The lower bound on $\meas{A_n},\ n=1,\ldots,N,$ 
which is to hold for any discretization implies that
\be \label{Nest}
N \e^3 \leq C;
\ee
as in all the estimates which follow $C$ denotes a positive constant
which may depend on $D$, $C^*$, or the restriction of $f$ on $[0,T]$.
We assume that the assumptions of 
Theorem~\ref{main} or Theorem~\ref{mainw} hold. 
In the situation
of Theorem~\ref{main} we redefine $\Ted$ such that 
$[0,\Ted]$ is the maximal interval on which
the estimate
\be \label{tedb}
\|\ovrho(t)\|_\infty \leq 2 D 
\ee
holds in addition to those stated in (\ref{teda}). 
In the arguments which follow
there is almost no need to distinguish between Theorem~\ref{main}
and Theorem~\ref{mainw}.
The only difference is the bound on $\ovrho$:
$\ovrho(t,r) \leq 2 D$ for $t \in [0,\Ted]$ in case of Theorem~\ref{main}
and $\ovrho(t,r) \leq C^\star$ for $t \in [0,\Ted]$ in case of 
Theorem~\ref{mainw}.

We will use the following abbreviation for the index set of the decomposition
of $\supp \fn$:
\[
\CN:=\{1,2,\ldots,N\}.
\]
We start our chain of auxiliary results by collecting some information
on our hat function $\chid$:
\begin{lemma} \label{chid}
Let $\d >0$ and define
\[ 
\chi(\zeta):=\frac{1}{\d}
\int_{-\infty}^\zeta \chid(\xi) \, ds,\ \zeta \in \R. 
\]
Then
\[ 
\left| \chi(\zeta) - \chi(\xi) \right| \leq \frac{1}{\d} |\zeta-\xi|,\
\ 0 \leq \chi(\zeta) \leq 1,\ \zeta,\ \xi \in \R,
\]
and for $a,b,\zeta \in \R$ with $|a-\zeta| \geq |a-b|+\d$,
\[
\chi(\zeta-a)=\chi(\zeta-b).
\]
\end{lemma}

The proof is straight forward.
Next we collect some additional information on the characteristic flow
of the Vlasov-Einstein system.
As before, we often abbreviate
$z=(x,v)\in \R^3\times \R^3$ and $Z=(X,V)$ accordingly.
\begin{lemma} \label{chinida}
For all characteristics starting in the support of $\fn$
and all $t \in [0,T]$, 
\[
|Z(t,0,z_1)-Z(t,0,z_2)| \leq C |z_1-z_2|,
\]
in particular, for $A \subset \supp \fn \subset \R^6$
and $A(t):=Z(t,0,A)$,
\[
\diam{A(t)} \leq C \diam{A}.
\]
Also
\beas
&&\bigl|(R,W)(t,0,r_1,w_1,L_1)-(R,W)(t,0,r_2,w_2,L_2)\bigr|\\
&&
\quad\quad\quad\quad\quad\quad\quad\quad\quad\quad
\leq C \bigl( |r_1-r_2| + |w_1-w_2| + | L_1-L_2| \bigr).
\eeas
\end{lemma}
\prf
The first estimate is included in \cite[Prop.~2.2]{RH}, 
with the bounds from Lemma~\ref{bounds} this implies the third estimate;
note that $1/R(t)$ is bounded.
\prfe

We now start establishing bounds on the approximating solutions;
recall that constants $C$ may depends on the true solution $f$,
but never on the discretization parameters $\e$ and 
$\d$.
\begin{lemma} \label{discrbounds}
For all $t \in [0,\Ted]$, $r > 0$ and $n \in \CN$,
\[
\ovrho(t,r),\; \ovp(t,r),\; |\ovj(t,r)|,\;\ovm(t,r), 
|\ovla(t,r)|,\; |\ovmu(t,r)|,\;
|\ovlap(t,r)|,\;  |\ovlas(t,r)|,\; |\ovmus(t,r)|
\leq C
\]
and
\[
\left|\dot{\ovM}_n(t) \right| \leq C \e^3.
\]
Moreover, for $r \not\in [1/(4D),3 D]$,
\[ 
\ovrho(t,r) = \ovp (t,r) = \ovj (t,r) = 0.
\]
\end{lemma}
\prf
Let $t \in [0,\Ted]$. By definition of $\Ted$,  
$r < 1/(4D)$ implies 
$r - \ovR_n(t) < 1/(4D) - 1/(2D) \leq - \d$,
and $r > 3 D$ implies $r - \ovR_n(t) > 3 D - 2 D \geq \d$. Thus
\be \label{rhob=0}
\ovrho(t,r) = 0,\ r \notin [1/(4D), 3 D].
\ee
A bound on $\ovrho$ of the desired type holds either by assumption
or by definition of $\Ted$.
It is obvious that 
$\ovp,\ |\ovj| \leq \ovrho$ so that the assertions for $\ovp$ and $\ovj$
follow.
Using the bounds on the approximating source terms and 
(\ref{rhob=0}) 
the remaining bounds follow from the definitions of the respective
quantities and of $\Ted$.
\prfe
 
Next we present an important tool for analyzing the 
approximating source terms:
\begin{lemma} \label{Icount}
Let $0 < r_0 < r_1$. For $t \in [0,\Ted]$ define
\[
I(t) :=
\left\{ n \in \N \mid R_n(t) \in [r_0,r_1] \right\},\
\ovI(t) 
:=
\left\{ n \in \N \mid \ovR_n(t) \in [r_0,r_1] \right\}.
\]
Then
\beas
|I(t)| 
&\leq& 
C \frac{r_1-r_0+\e}{\e^3},\\
|\ovI(t)| 
&\leq& 
C \frac{r_1-r_0+\NormRt+\e}{\e^3},
\eeas
where $|I|$ denotes the number of elements of a set
$I \subset \N$. Moreover,
\[
\sum_{n \in \ovI(t)} \ovM_n(t) \leq C ( r_1 - r_0 + \d ).
\]
\end{lemma}
\prf
The proof of the first assertion follows the one of \cite[Lemma~2]{GS}.
Fix $n \in I(t)$ and $(r,w,L) \in A_n$.
Then Lemma \ref{chinida} implies that  
\[
R(t,0,r,w,L)| \in [r_0-C\e,r_1+C\e],
\]
and thus
\[ 
\bigcup_{n \in I(t)} A_n(t) \subset [r_0-C \e, r_1 + C \e] \times [-D,D]
\times [0,D] .
\]
Using Lemma~\ref{phspvol}, the bound on $\l$
from Lemma~\ref{bounds}, and the lower bound on $\meas{A_n}$ we find
\[
(r_1-r_0+2 C \e) 2 D^2 
\geq
\sum_{n \in I(t)} \meas{A_n(t)} 
\geq \frac{1}{C}\sum_{n \in I(t)} \meas{A_n} 
\geq C |I(t)| \e^3,
\]
which implies the first assertion. 
The second one follows from the first, since
\[ 
\ovI(t) \subset \left\{ n \in \CN \mid 
R_n(t) \in [r_0-\NormRt,r_1+\NormRt] \right\}.
\]
To prove the third inequality define for $r \geq \d/2$,
\[ 
\ovI(t,r) := \left\{ n \in \CN \mid 
\ovR_n(t) \in [r-\d/2, r+\d/2] \right\}.
\]
Then $n \in \ovI(t,r)$ implies $|\ovR_n(t)-r| \leq \d/2$ and 
therefore $\chid(\ovR_n(t)-r) \geq 1/2$. Hence
\beas
\sum_{n \in \ovI(t,r)} \ovM_n(t) 
&\leq& 
2 \sum_{n \in \ovI(t,r)} \ovM_n(t)
        \chid(\ovR_n(t)-r)\\        
&\leq& 
8 \pi r^2 \d \, \ovrho(t,r) 
\leq 8 \pi r^2\,\d \max\{C^\star, 2 D \} 
= C r^2 \d,\ r \geq \d/2.
\eeas
Now define 
$k:=\left\lceil (r_1-r_0)/\d\right\rceil$ and $s_i := r_0 + i\d$,
$i=0, \dots,k$. Here $\lceil \zeta \rceil$ denotes the smallest integer larger
than or equal to $\zeta$. Then $s_0 = r_0$ and $s_k \geq r_1$, and therefore
\[ 
\ovI(t) \subset \bigcup_{i=0}^{k-1} \ovI(t,s_i+\d/2).
\]
Using this inclusion and the previous estimate for every $i$ 
together with the fact that by definition of $\Ted$ only $r_1\leq 3D$
needs to be considered we get the estimate 
\beas
\sum_{n \in \ovI(t)} \ovM_n(t) 
&\leq& 
\sum_{i=0}^{k-1} 
\sum_{n \in \ovI(t,s_i+\d/2)} \ovM_n(t)
\leq 
\sum_{i=0}^{k-1} C (s_i+\d/2)^2 \,\d \\
&\leq&
C (r_1+\d/2)^2 \, k \,\d +
\leq 
C \left\lceil (r_1-r_0)/\d\right\rceil \d 
\leq C ( r_1 - r_0 + \d ),
\eeas
and the proof is complete.
\prfe

In order to analyze the errors appearing at the
source terms and the metric coefficients due to the discretization we
define quantities that are intermediate between the true quantities
and their approximations in the sense that the formulas for the latter source
terms are evaluated at the true characteristics with the true weight functions.
These intermediate
quantities are denoted by double-bars:
\[ 
\ovvrho(t,r) := \frac{1}{4 \pi r^2 \d} \sum_{n=1}^N 
E_n(t) M_n(t) \chid(r-R_n(t))
\]
where analogously to (\ref{discredef}),
\[
E_n(t):=\sqrt{1+W_n^2(t) + L_n/R_n^2(t)},
\]
with corresponding definitions for $\ovvp$ and $\ovvj$.
Moreover, $\ovvM_n$ is to be the solution of
\[ 
\dot{\ovvM}_n(t)= -\left(\dot \lambda(t,R_n(t))+\dot{R}_n(t)
\lambda'(t,R_n(t)) \right) \ovvM_n(t), \quad \ovvM_n(0) = M_n(0).
\]
Note that we have, similar to $\ovrho$ and $\ovM$, $\|\ovvrho(t)\|_\infty +
\e^{-3} \ovvM(t) \leq C$ for all $t \in [0,T]$ and $\ovvrho(t,r) =0$ for all
$r \not\in [1/(4D),3 D]$, $t \in [0,T]$.
\begin{lemma} \label{mu-mub}
For all $t \in [0,\Ted]$,
\beas
\|\mu(t)-\ovmu(t)\|_\infty 
&+&
\|\l(t)-\ovla(t)\|_\infty + \|m(t)-\ovm(t)\|_\infty \\
&\leq& 
C \Bigl(\|\r(t)-\ovvrho(t)\|_\infty +
         \|p(t)-\ovvp(t)\|_\infty \\
&&\qquad {} + \NormRt + \NormWt + \NormMt \Bigr).
\eeas
\end{lemma}
\prf
Let $t \in [0,\Ted]$ and $r > 0$. We start with analyzing the error 
$m-\ovm$. 
Due to the definition of $\ovm$, $\ovrho$, and $\ovvrho$,
\beas
&&
m(t,r) - \ovm(t,r) 
= 
4 \pi \int_0^r (\r(t,s)-\ovrho(t,s))  s^2 ds
=
4 \pi \int_0^r (\r(t,s)-\ovvrho(t,s)) s^2 ds\\
&&
\qquad {}+  \frac{1}{\d} \int_0^r \sum_{n=1}^N \Bigl( E_n(t)
         M_n(t) \chid(s-R_n(t)) - \ovE_n(t) \ovM_n(t)   
        \chid(s-\ovR_n(t)) \Bigr) ds \\
&&
\quad =: F_1 + F_2.
\eeas
Since for $r \geq 3 D$ we have $\r(t,r)=\ovvrho(t,r)=0$,
\[
|F_1| \leq C \|\r(t)-\ovvrho(t)\|_\infty .
\]
Using the definition of $\chi$, cf.\ Lemma~\ref{chid}, we have
\beas
|F_2| 
&\leq&
\left| \sum_{n=1}^N ( M_n(t) - \ovM_n(t) )
E_n(t) \chi(r-R_n(t)) \right| \\
&&{}
+\left| \sum_{n=1}^N \ovM_n(t) \left(E_n(t) -
        \ovE_n(t) \right) \chi(r-R_n(t)) \right| \\
&&{}
+\left| \sum_{n=1}^N \ovM_n(t) \ovE_n(t) 
\left( \chi(r-R_n(t)) - \chi(r-\ovR_n(t)) \right) \right| =: 
F_{21} + F_{22} + F_{23}.
\eeas
Since $E_n$ and $\chi$ are bounded and $\e^3 N \leq C$ by (\ref{Nest}),
\[
F_{21} \leq C \sum_{n=1}^N |M_n(t) - \ovM_n(t)| \leq C \NormMt.
\]
Since $W_n(t)$, $\ovW_n(t)$, $1/R_n(t)$, $1/\ovR_n(t)$,
and $L_n$ are bounded,
\[
\left|E_n(t)-\ovE_n(t)\right|
\leq C \left( |R_n(t)-\ovR_n(t)| +|W_n(t)-\ovW_n(t)| \right),
\]
and thus 
\beas
F_{22} 
&\leq& 
C \sum_{n=1}^N \ovM_n(t) \bigl(\NormRt + \NormWt\bigr) \\
&\leq& 
C \bigl(\NormRt + \NormWt\bigr).
\eeas
If we define
\[ 
I_1(t):= 
\left\{ n \in \CN \mid |R_n(t)-r| \leq \|R(t) - \ovR(t)\|_\infty + \d \right\}
\]
then by Lemma~\ref{chid},
\[ 
\chi(r-R_n(t)) = \chi(r-\ovR_n(t)),\  n\not\in I_1(t),
\]
and by Lemma \ref{Icount},
\[ 
|I_1(t)| \leq C \frac{\|R(t) - \ovR(t)\|_\infty + \d}{\e^3}
\]
which implies, again with Lemma~\ref{chid}, that
\[
F_{23} 
\leq 
C \sum_{n \in I_1(t)} \ovM_n(t) \min\left\{ 
        \frac{|R_n(t)-\ovR_n(t)|}{\d},1\right\}\\
\leq C \NormRt.
\]
This completes the desired estimate for $m-\ovm$. 
The estimate for $\l-\ovla$ follows directly from the expressions for
$\l$ and $\ovla$ in terms of $m$ and $\ovm$, the fact that 
$\l(t,r) = \ovla(t,r)=0$ for $r<1/(4D)$, and the estimate which we just
completed. 

It remains to consider $\m - \ovmu$.
For $r < 1/(4D)$ we have $\mu'(t,r)=\ovmus(t,r)=0$ and therefore
$|\mu(t,r)-\ovmu(t,r)| = |\mu(t, 1/(4D))-\ovmu(t, 1/(4D))|$. 
Let $r \geq  1/(4D)$. Then
the definitions of $\mu$, $\ovmu$, $\mu'$, and $\ovmu'$ yield
\beas
|\mu(t,r)-\ovmu(t,r)| 
&=&
\left|\int_r^\infty \ovmus(t,s)-\mu'(t,s)\, ds\right|\\
&\leq& 
\int_r^\infty e^{2 \ovla(t,s)}
\left|\frac{\ovm(t,s)}{s^2} - \frac{m(t,s)}{s^2} \right| ds \\
&&
{}+ \int_r^\infty \left|e^{2 \ovla(t,s)}-e^{2 \l(t,s)}\right| 
\left( \frac{m(t,s)}{s^2} + 4 \pi s p(t,s) \right) ds \\
&&
{}+ \left|\int_r^\infty  4 \pi s e^{2 \ovla(t,s)}
\bigl(\ovp(t,s) - p(t,s)\bigr) ds \right|\\
&\leq& 
C \|m(t)-\ovm(t)\|_\infty + C \|\l(t)-\ovla(t)\|_\infty \\
&&
{}+ \left|\int_r^\infty 
4 \pi s e^{2 \ovla(t,s)} \bigl(\ovp(t,s) - p(t,s) \bigr) ds \right|.
\eeas
We split the difference
$\ovp-p$ into $\ovp - \ovvp$ and $\ovvp-p$. The integral resulting 
from the second term can obviously be estimated as desired.
As to the first term,
\beas
&&
\left| \int_r^\infty 4 \pi s e^{2 \ovla(t,s)} ( \ovvp(t,s) - \ovp(t,s) ) 
ds\right|\\
&&\qquad \leq
\frac{C}{\d} \sum_{n=1}^N |M_n(t) - \ovM_n(t)|
\int_{1/(4D)}^{3D} \chid(s-R_n(t))\,ds\\
&&
\qquad \quad{}+ \frac{C}{\d} \sum_{n=1}^N \ovM_n(t) 
\left|\frac{W_n^2(t)}{E_n(t)} - \frac{\ovW_n^2(t)}{\ovE_n(t)}\right|
\int_{1/(4D)}^{3D} \chid(s-R_n(t)) ds\\
&&
\qquad \quad{}+ \frac{C}{\d} 
\sum_{n=1}^N \ovM_n(t) 
\left|\int_r^\infty \frac{e^{2 \ovla(t,s)}}{s} 
\left( \chid(s-R_n(t)) - \chid(s-\ovR_n(t)) \right) ds\right|\\
&&
\qquad =:
F_3 + F_4 + F_5
\eeas
where we used the bounds on the supports of $\ovvp$ and $\ovp$ and recalled
that $r\geq 1/(4D)$.
The first two terms are comparatively easy to deal with:
Define $I_2(t,s) := \left\{ n \in \CN \mid |s-R_n(t)|\leq \d\right\}$. 
By Lemma~\ref{Icount},
$|I_2(t,s)| \leq C (\d+\e)\e^{-3} \leq C \d \e^{-3}$
for all $s \geq r$. Thus
\beas
F_3 
&\leq& 
\frac{C}{\d} \e^3 \NormMt \int_{1/(4D)}^{3 D}   |I_2(t,s)|ds
\leq 
C \NormMt,\\
F_4 
&\leq& 
\frac{C}{\d} \e^3 \left( \NormRt +\NormWt \right)
\int_{1/(4D)}^{3 D}   |I_2(t,s)| ds\\
&\leq& 
C \left( \NormRt +\NormWt \right).
\eeas
To deal with $F_5$ 
we fix some $n \in \CN$ and assume without loss of generality
that $\ovR_n(t) \leq R_n(t)$. A change of variables yields
\beas
&&
\left|\int_r^\infty \frac{e^{2 \ovla(t,s)}}{s}
\left(\chid(s-R_n(t)) -\chid(s-\ovR_n(t))\right) ds\right|\\
&&
\qquad
=\left|\int_{r-R_n(t)}^\infty 
\frac{e^{2 \ovla(t,z+R_n(t))}}{z+R_n(t)} \chid(z)\, dz 
-\int_{r-\ovR_n(t)}^\infty \frac{e^{2 \ovla(t,z+\ovR_n(t))}}{z+\ovR_n(t)}
\chid(z)\, dz\right|\\
&&
\qquad
\leq\int_{r-R_n(t)}^{r-\ovR_n(t)}\frac{e^{2 \ovla(t,z+R_n(t))}}{z+R_n(t)}  
\chid(z)\, dz \\
&&
\qquad \quad
{} +\int_{r-\ovR_n(t)}^\infty 
\left|\frac{e^{2 \ovla(t,z+R_n(t))}}{z+R_n(t)} - 
\frac{e^{2 \ovla(t,z+\ovR_n(t))}}{z+\ovR_n(t)}\right|
\chid(z) \, dz\\
&&\qquad
\leq
C \int_{r-R_n(t)}^{r-\ovR_n(t)}\chid(z)\, dz
+ C \int_{r-\ovR_n(t)}^\infty \chid(z) \, dz \left|R_n(t) - \ovR_n(t)\right|\\
&&\qquad
\leq
C \min\left\{ |R_n(t)-\ovR_n(t)|,\d \right\} + C \d \NormRt.
\eeas
Define
\beas
I_3(t) 
&:=& 
\left\{n \in \CN \mid [r-R_n(t),r-\ovR_n(t)] \cap 
[-\d,+\d] \not = \emptyset\right\}\\
&\subset& 
\left\{ n \in \CN \mid |R_n(t)-r| \leq \d + \NormRt \right\}.
\eeas
Then Lemma~\ref{Icount} yields $|I_3(t)| \leq C \e^{-3} (\d + \NormRt )$,
and putting everything together, we obtain
\beas
|F_5| 
&\leq& 
\frac{C \e^3}{\d} 
\sum_{n \in I_3(t)} \min\left\{ |R_n(t)-\ovR_n(t)|,\d \right\}
+ \frac{C \e^3}{\d} \sum_{n=1}^N \d \NormRt\\
&\leq& 
C \NormRt
\eeas
by (\ref{Nest}), and the proof is complete.
\prfe

Next we estimate the differences of the intermediate, double-barred source
terms and their true counterparts:
\begin{lemma} \label{rho-rhobb}
For all $t \in [0,\Ted]$,
\beas 
\|\r(t)-\ovvrho(t)\|_\infty
&\leq& 
C \left( \d + \e/\d \right),\\
\|\r(t)-\ovvrho(t)\|_\infty
&\leq& 
C \left( \d^2 + 
\e/\d \right), \ \mbox{provided}\ \r'' \in C([0,T]\times\R^+),
\eeas
and the analogous assertions hold for $p$ and $\j$.
\end{lemma}
\prf
We restrict ourselves to considering $\r$;
$p$ and $\j$ can be dealt with in the same fashion.
Let $t \in [0,\Ted]$ and $r>0$. Then 
\beas
4 \pi r^2 \d\, \r(t,r) 
&=& 
\int_{r-\d}^{r+\d} 4 \pi s^2 \r(t,s) \chid(r-s)\, ds\\
&&
{}+ \int_{r-\d}^{r+\d} 
4 \pi \left( r^2 \r(t,r) -
s^2 \r(t,s) \right) \chid(r-s)\, ds \\
&=:&
F_1 + F_2,\\
4 \pi r^2 \d\, \ovvrho(t,r) 
&=& 
\sum_{n=1}^N E_n(t) \chid(r-R_n(t)) \int_{A_n(t)} f(t,z) \, dz\\
&=& 
\sum_{n=1}^N \int_{A_n(t)} E f(t,z) \chid(r-|x|)\, dz\\
&& 
{}+\sum_{n=1}^N \int_{A_n(t)} \left( E_n(t) -E \right) 
f(t,z) \chid(r-R_n(t))\, dz \\
&&
{} +\sum_{n=1}^N \int_{A_n(t)} E f(t,z) 
\left( \chid(r-R_n(t))-\chid(r-|x|) \right) \, dz \\
&=:& 
F_3 + F_4 + F_5.
\eeas
Clearly, $F_1 = F_3$. Since $\r(t,r) = \ovvrho(t,r) = 0$ for all 
$r \not\in [1/(4D), 3 D]$, 
\[ 
\left|\r(t,r)-\ovvrho(t,r)\right| 
\leq \frac{C}{\d} 
\bigl( |F_2| + |F_4| + |F_5| \bigr).
\]
Using the fact that we assume a bound on
the derivative of $s^2 \r(t,s)$ with respect
to $s$ we find 
\[
|F_2| \leq C \int_{r-\d}^{r+\d} |r-s| \chid(r-s)\, ds 
\leq C \d^2.
\]
To estimate $F_4$ let $z=(x,v) \in A_n(t)$. Let $Z_n(t)=(X_n(t),V_n(t))$
be any Cartesian characteristic starting at a point
with $(r,w,L)$ coordinates $(R_n(0),W_n(0),L_n)$. Then $Z_n(t) \in A_n(t)$
and $|z-Z_n(t)| \leq \diam{A_n(t)} \leq C \e$ by Lemma~\ref{chinida}.
The function $E$ used to compute $\r$ and $\ovvrho$ is Lipschitz
on the relevant domain; the same is true for $p$ and $\j$.  
If we define 
$I_1(t):=\left\{ n \in \CN \mid |R_n(t) -r| \leq \d\right\}$ then
Lemma~\ref{Icount} yields $|I_1(t)| \leq C \d \e^{-3}$ and therefore,
\beas
|F_4| 
&\leq& 
\sum_{n\in I_1(t)} \int_{A_n(t)} \left| E_n(t)-E \right| 
f(t,z) \, dz\\
&\leq& 
\sum_{n\in I_1(t)} C \diam{A_n(t)} \int_{A_n(t)} f(t,z) \, dz
\leq C \d \e.
\eeas
For $F_5$ we define 
\[
I_2(t):=\left\{ n \in \CN \mid 
|R_n(t)-r|\leq \d+\max_{m\in \CN}\diam{A_m(t)} \right\}.
\] 
Then $\chid(r-R_n(t))=\chid(r-|x|)$ for all $n \not\in
I_2(t)$ and $z\in A_n(t)$, and by Lemma~\ref{Icount}, 
$|I_2(t)| \leq C \e^{-3} (\d +\max_{m \in \CN} 
\diam{A_m(t)} + \e ) \leq C \e^{-3} \d$. Hence
\beas
|F_5| 
&\leq& 
\frac{1}{\d} \sum_{n\in I_2(t)} \int_
{A_n(t)} E f(t,z) |R_n(t)-|x|| dz\\
&\leq& 
\frac{C}{\d} \sum_{n\in I_2(t)} \diam{A_n(t)} M_n(t) \leq C \e,
\eeas
and the first assertion is proven. 

To prove the second one, we have to 
sharpen the estimate for $F_2$ in the case that $\r''$ exists and is 
continuous and thus bounded on the relevant domain.
 First we rewrite this term;
since throughout the argument $t$ remains
fixed we suppress this variable:
\beas
F_2 
&=& 
4 \pi \int_{r-\d}^{r+\d} (r^2-s^2) \r(r) \chid(r-s) \, ds\\
&&
{}+4 \pi \int_{r-\d}^{r+\d} 
\left( \r(r)-\frac{1}{2}\Bigl(\r(r+\d)+\r(r-\d)\Bigr)
\right) s^2 \chid(r-s) \, ds\\
&&
{}+4 \pi \int_{r-\d}^{r+\d} 
\left( \frac{1}{2}\Bigl(\r(r+\d)+\r(r-\d)\Bigr)-\r(s)\right) 
s^2 \chid(r-s) \, ds\\
&=:& 
F_6 + F_7 + F_8.
\eeas
By the mean value theorem
\[ 
|F_6| \leq  C \left|\int_{r-\d}^{r+\d} (r^2-s^2) \, ds\right| = C \d^3.
\]
To analyze $F_7$ we use the extra regularity of $\r$. By Taylor 
expansion,
\[
|F_7| \leq C\int_{r-\d}^{r+\d}\left| \r(r) - 
\frac{1}{2} \Bigl( \r(r+\d)+\r(r-\d)\Bigr) \right|s^2 \, ds 
\leq C \d^3.
\]
Finally, by the mean value theorem and Taylor expansion,
\beas
|F_8| 
&\leq& 
C\left|\int_{r-\d}^{r+\d} \left( \frac{1}{2} 
\Bigl( \r(r+\d)+\r(r-\d) \Bigr) - \r(s) \right)\, ds\right| \\
&=&
C \left|\d \r(r+\d) + \d \r(r-\d) - \int_r^{r+\d} \r(s)\, ds 
- \int_{r-\d}^r \r(s)\, ds \right| \\      
&\leq& 
C \left|\d \r(r+\d) + \d \r(r-\d) - 
\int_r^{r+\d}(\r(r+\d) + \r'(r+\d)(s-r-\d))\,ds \right.\\
&&
{}\qquad \qquad \qquad \qquad \qquad \left. -
\int_{r-\d}^r(\r(r-\d) + \r'(r-\d)(s-r+\d))\,ds \right| \\      
&&
{} + C \int_r^{r+\d} (s-r-\d)^2 ds + C \int_{r-\d}^r (s-r+\d)^2 ds\\
&=&
C \left|\frac{1}{2} \d^2 \r'(r+\d) - \frac{1}{2} \d^2 \r'(r-\d)\right|
+ C \d^3\\
&\leq&
C \d^3,
\eeas
and the proof is complete. \prfe

The next lemma will complement the previous one in that we
now estimate the differences between the approximations of the source
terms and the corresponding intermediate double-barred quantities.
It is at this point that we loose one order of $\d$ in the error 
estimates for the sources compared to the other quantities 
in Theorem~\ref{main}:
\begin{lemma} \label{rhob-rhobb}
For all $t \in [0,\Ted]$,
\[
\|\ovrho(t)-\ovvrho(t)\|_\infty
\leq C \left( \frac{\NormRt}{\d} + \NormWt + \NormMt \right),
\]
in particular, $\|\ovrho(0)-\ovvrho(0)\|_\infty= 0$. The same estimates
hold for $p$ and $\j$.
\end{lemma}
\prf
Since $\ovrho(t,r)=\ovvrho(t,r)=0$ for 
$r \not\in [1/(4D), 3 D]$ we only need to consider
$r \in [1/(4D), 3 D]$. By definition
\beas
\ovrho(t,r)-\ovvrho(t,r) 
&=&
\frac{1}{4 \pi r^2 \d} \sum_{n=1}^N \ovE_n(t) \ovM_n(t)
\left( \chid(r-\ovR_n(t)) - \chid(r-R_n(t)) \right) \\
&&
{}+ \frac{1}{4 \pi r^2 \d} \sum_{n=1}^N 
\left(\ovE_n(t) \ovM_n(t) - E_n(t) M_n(t) \right)\chid(r-R_n(t)).
\eeas
Let
$I_1(t) := \left\{ n \in \CN \mid |R_n(t)-r| \leq \d \right\}$ 
and
$I_2(t) := \left\{ n \in \CN \mid |\ovR_n(t)-r| \leq \d \right\}$.
Then Lemma~\ref{Icount} implies that $|I_1(t)| \leq C \e^{-3} \d$
and $\sum_{n\in I_2(t)} \ovM_n(t) \leq C \d$. Moreover,
$n \not\in I_1(t) \cup I_2(t)$ implies 
$\chid(r-\ovR_n(t)) = 0 = \chid(r-R_n(t))$. Using similar estimates
as in the lemmas above we find
\beas
&&
|\ovrho(t,r)-\ovvrho(t,r)|
\leq \frac{C}{\d} \sum_{n\in I_1(t) \cup I_2(t)}
\ovM_n(t)\; |R_n(t)-\ovR_n(t)|/\d \\
&&
\quad +\frac{C \e^3}{\d} \!\!\sum_{n\in I_1(t)} 
\left(|R_n(t)-\ovR_n(t)|+|W_n(t)-\ovW_n(t)|+\e^{-3} |M_n(t)-\ovM_n(t)|\right)\\
&&
\leq C \frac{\NormRt}{\d^2} \left( \sum_{n\in I_1(t)} \e^3
+ \sum_{n\in I_2(t)} \ovM_n(t) \right)\\
&&
\quad + \frac{C \e^3}{\d} | I_1(t)| \left( \NormRt
+\NormWt + \NormMt \right),
\eeas
and the proof is complete.
\prfe

We now turn to the analysis of the evolution equations.
\begin{lemma} \label{m-mbb}
For all $t\in [0,\Ted]$,
\[
\NormMqqt := \e^{-3} \max_{n \in \CN} 
\sup_{s \in [0,t]} |M_n(s)-\ovvM_n(s)|
\leq C \e. 
\]
\end{lemma}
\prf
Let $t\in [0,\Ted]$ and $n \in \CN$. The definition of $\ovvM_n$ and
Lemma~\ref{phspvol} with $A=A_n$ yields
\beas
M_n(t)-\ovvM_n(t) 
&=& 
\int_0^t \left( \dot\l(s,R_n(s))
+ \dot R_n(s) \l'(s,R_n(s)) \right) \ovvM_n(s) \,ds\\
&&
{}- \int_0^t \int_{A_n(s)}
f(s,z) \left( \dot\l(s,|x|) + \dot R(s,0,Z(0,s,z))\l'(s,|x|)\right) 
\, dz\, ds\\
&=& 
\int_0^t \left[ \left( \dot\l(s,R_n(s))
+ \dot R_n(s) \l'(s,R_n(s)) \right) \left(\ovvM_n(s)-M_n(s)\right)\right]\,ds\\
&&
{}+\int_0^t \int_{A_n(s)}
f(s,z)\Bigl[ \dot\l(s,R_n(s))-\dot\l(s,|x|)+\dot R_n(s)\l'(s,R_n(s))\\
&&
\qquad\qquad\qquad\qquad\quad\quad\quad\quad
{}-\dot R(s,0,Z(0,s,z))\l'(s,|x|)\Bigr] \, dz\, ds.
\eeas
For $z=(x,v) \in A_n(s)$,
\[
|\dot\l(s,R_n(s))-\dot\l(s,|x|)| \leq C \bigl| R_n(s)-|x| \bigr| 
\leq C \diam{A_n(s)} \leq C \e,
\]
and
\beas
\bigl|\dot R_n(s)\l'(s,R_n(s))
&-&
\dot R(s,0,Z(0,s,z))\l'(s,|x|)\bigr|\\
&\leq& 
\left|\dot R(s,0,z_n)-\dot R(s,0,Z(0,s,z))\right| \left|\l'(s,R_n(s))\right| \\
&&
{}+ \left|\dot R(s,0,Z(0,s,z))\right| \bigl|\l'(s,R_n(s))-\l'(s,|x|)\bigr|\\
&\leq& 
C \bigl( |z_n-Z(0,s,z)| + |R_n(s)-|x|| \bigr)\\
&\leq& 
C \bigl( \diam{A_n(0)} + \diam{A_n(s)} \bigr) \leq C \e.
\eeas
Hence,
\beas
\left|M_n(t)-\ovvM_n(t)\right| 
&\leq& 
C \int_0^t \left|M_n(s)-\ovvM_n(s)\right| \, ds + C \e \int_0^t M_n(s) \, ds\\
&\leq& 
C \int_0^t \left|M_n(s)-\ovvM_n(s) \right| \, ds + C \e^4,
\eeas
and a Gronwall argument completes the proof.\prfe

\begin{lemma} \label{dr-drb}
For all $t\in [0,\Ted]$ and $n \in \CN$,
\beas 
|\dot R_n(t) - \dot{\ovR}_n(t)| 
&\leq& C \Big(|R_n(t)-\ovR_n(t)| + |W_n(t)-\ovW_n(t)| \\
&&
\qquad + \|\mu(t)-\ovmu(t)\|_\infty + \|\l(t)-\ovla(t)\|_\infty \Big).
\eeas
\end{lemma}
\prf
By the equations for $\dot R_n(t)$ and $\dot{\ovR}_n(t)$ and straight forward
estimates,
\beas
|\dot R_n(t) - \dot{\ovR}_n(t)|
&\leq& 
\left| e^{(\mu-\l)(t,R_n(t))} - e^{(\mu-\l)(t,\ovR_n(t))}\right| \\
&&
{}
+\left|e^{(\mu-\l)(t,\ovR_n(t))}-e^{(\ovmu-\ovla)(t,\ovR_n(t))}\right|
+ \left| \frac{W_n(t)}{E_n(t)} - \frac{\ovW_n(t)}{\ovE_n(t)} \right|,
\eeas
and the assertion follows.
\prfe

The use of the following technical lemma will become obvious later:
\begin{lemma} \label{gmn}
Let $g_{n,m} \in C^1([0,\Ted]),\ m,n \in \CN$. Then for
all $t \in [0,\Ted]$ and $m\in \CN$,
\beas
&&
\Bigg| \frac{1}{\d} \sum_{n=1}^N \int_0^t g_{n,m}(s) 
\Big[
\chid(\ovR_m(s)-R_n(s)) \left( \dot{\ovR}_m(s)-\dot R_n(s) \right)\\
&&
\quad\quad\quad\quad\quad\quad\quad\quad\quad\quad\quad\quad\quad
-\chid(\ovR_m(s)-\ovR_n(s))\left( \dot{\ovR}_m(s)-\dot{\ovR}_n(s)\right)
\Big] \, ds \Bigg|\\
&&
\leq C \, \e^{-3} \max_{k,l \in \CN}
\bigl( \|g_{k,l}\|_\infty + \|\dot g_{k,l}\|_\infty\bigr)\,\NormRt.
\eeas
(Note that $C$ does not depend on the functions $g_{m,n}$.)
\end{lemma}
\prf
Let $m \in \CN$ be fixed, denote the left hand side of the inequality by $S$, 
and define
\[ 
d_n(s) := \ovR_m(s)-R_n(s),\ \ovd_n(s) := \ovR_m(s)-\ovR_n(s),\ 
s \in [0,\Ted].
\]
Integration by parts and $d_n(0)=\ovd_n(0)$ implies that
\[
S
= \sum_{n=1}^N \left[ g_{n,m}(t) \left( \chi(d_n(t)) - 
\chi(\ovd_n(t)) \right) - \!\!\int_0^t \dot g_{n,m}(s) 
\left(\chi(d_n(s)) - \chi(\ovd_n(s)) \right) \, ds \right].
\]
For $s \in [0,\Ted]$ define
\[ 
I(s):=\left\{ n \in \CN \mid 
|R_n(s)- \ovR_m(s)| \leq \d +\NormRs \right\}.
\]
By Lemma \ref{Icount}, 
\[
|I(s)| \leq C \e^{-3} \bigl(\d + \NormRs + \e\bigr)
\leq C \e^{-3} \bigl(\d + \NormRt \bigr),\
s \in [0,\Ted].
\] 
If $R_n(s) < \ovR_m(s)-\d-\NormRs$ then $d_n(s)\geq \d$ and
$\ovd_n(s) \geq \d$, and hence
$\chi(d_n(s)) = 1 = \chi(\ovd_n(s))$.
If $R_n(s) > \ovR_m(s)+\d+\NormRs$ then $d_n(s)\leq - \d$ and
$\ovd_n(s) \leq -\d$, and hence
$\chi(d_n(s)) = 0 = \chi(\ovd_n(s))$.
Hence, using Lemma~\ref{chid},
\beas
|S| 
&\leq& 
\sum_{n \in I(t)} |g_{n,m}(t)| \left| \chi(d_n(t)) - \chi(\ovd_n(t)) \right|\\
&&
{}+ \int_0^t \max_{k,l \in \CN}
|\dot g_{k,l}(s)| \sum_{n \in I(s)} \left| \chi(d_n(s)) - 
\chi(\ovd_n(s)) \right| \, ds\\
&\leq& 
\max_{k,l \in \CN} \|g_{k,l}\|_\infty
|I(t)| \min\left\{\NormRt/\d,1\right\}\\
&&
{} + \max_{k,l\in\CN} \|\dot g_{k,l}\|_\infty
\int_0^t |I(s)| \min\left\{\NormRs/\d,1\right\} \, ds\\
&\leq& 
C \max_{k,l \in \CN}
\bigl( \|g_{k,l}\|_\infty + \|\dot g_{k,l}\|_\infty \bigr)\,
\bigl(\d+\NormRt\bigr)\,\e^{-3}\\
&&
{}  \qquad \min\left\{\NormRt/\d,1\right\},
\eeas
and the assertion follows.
\prfe

\begin{lemma} \label{mbb-mb}
Let 
\beas
e(s)
&:=&
\|\r(s)-\ovvrho(s)\|_\infty + \|p(s)-\ovvp(s)\|_\infty + 
\|\j(s)-\ovvj(s)\|_\infty\\
&&
{}+ \NormWs + \NormMs.
\eeas 
Then
\[
\|\ovvM(t)-\ovM(t)\| \leq C \int_0^t e(s) \, ds + C \, \NormRt,\ 
t \in [0,\Ted].
\]
\end{lemma}
\prf
Let $t \in [0,\Ted]$ and $m \in \CN$. By the definitions of $\ovM_m$
and $\ovvM_m$,
\beas
&&
\ovvM_m(t) - \ovM_m(t)\\ 
&&\qquad
= 
\int_0^t \Bigl[ \ovlap(s,\ovR_m(s)) + 
\dot{\ovR}_m(s) \ovlas(s,\ovR_m(s)) \Bigr] \ovM_m(s) \, ds\\
&&\qquad\quad
{} - \int_0^t \Bigl[ \dot\l(s,R_m(s)) + 
\dot{R}_m(s) \l'(s,R_m(s)) \Bigr] \ovvM_m(s) \, ds\\
&&\qquad
=
\int_0^t \Bigl[ \ovlap(s,\ovR_m(s)) + 
\dot{\ovR}_m(s) \ovlas(s,\ovR_m(s))\\
&&\qquad\quad\quad\quad\quad - \dot\l(s,\ovR_m(s)) - 
\dot{\ovR}_m(s) \l'(s,\ovR_m(s)) \Bigr] \, \ovM_m(s) \, ds\\
&&\qquad\quad
{}+\int_0^t \Bigl[\dot\l(s,\ovR_m(s)) - \dot\l(s,R_m(s))
+\left( \dot{\ovR}_m(s) - \dot R_m(s) \right) \l'(s,\ovR_m(s))\\
&&\qquad\quad\quad\quad\quad 
+\dot{R}_m(s) \left( \l'(s,\ovR_m(s))-\l'(s,R_m(s)) \right)  \Bigr] \, 
\ovM_m(s) \, ds\\
&&\qquad\quad
{}+\int_0^t \Bigl[ \dot\l(s,R_m(s)) + 
\dot{R}_m(s) \l'(s,R_m(s)) \Bigr] \left( \ovM_m(s)-\ovvM_m(s)
\right) \, ds\\
&&\qquad
=: 
F_1 + F_2 + F_3.
\eeas
By Lemma~\ref{bounds}, 
\[
|F_3| \leq C \int_0^t \bigl|\ovM_m(s) - \ovvM_m(s)\bigr|\, ds.
\]
By Lemma~\ref{bounds} and Lemma~\ref{dr-drb},
\beas
|F_2| 
&\leq& 
C \e^3 \int_0^t 
\bigl(|R_m(s)-\ovR_m(s)| + |\dot R_m(s)-\dot{\ovR}_m(s)|\bigr)\, ds\\
&\leq& 
C \e^3 \int_0^t \Bigl(\|\m(s)-\ovmu(s)\|_\infty
+ \|\l(s)-\ovla(s)\|_\infty\\
&&
\quad\quad\quad\quad\quad+ \NormRs + \NormWs \Bigr) \, ds.
\eeas
Now we use the formulas for the derivatives of $\l$, $\ovla$
respectively where we smuggle in the intermediate
double-barred source terms to estimate $F_1$:
\beas
|F_1| 
&=& 
\Bigg| 4 \pi \int_0^t \ovR_m(s) \ovM_m(s)\\ 
&&
\qquad \qquad
\Biggl[e^{(\l+\m)(s,\ovR_m(s))} 
\j(s,\ovR_m(s)) - e^{(\ovla+\ovmu)(s,\ovR_m(s))} \ovj(s,\ovR_m(s))\\
&&
\qquad\qquad\quad+ \dot{\ovR}_m(s) 
\left(e^{2\ovla(s,\ovR_m(s))}
\left(\ovrho(s,\ovR_m(s))- \frac{\ovm(s,\ovR_m(s))}{\ovR_m^3(s)}\right)
\right.\\
&&
\qquad\qquad\qquad\qquad\qquad \left. - e^{2\l(s,\ovR_m(s))}
\left(\r(s,\ovR_m(s))- \frac{m(s,\ovR_m(s))}{\ovR_m^3(s)}\right)
\right)\Biggr]  \, ds \Bigg|\\
&\leq& 
C \e^{-3}
\int_0^t \left| 
e^{2\ovla(s,\ovR_m(s))}\frac{\ovm(s,\ovR_m(s))}{\ovR_m^3(s)}
-  e^{2\l(s,\ovR_m(s))}\frac{m(s,\ovR_m(s))}{\ovR_m^3(s)}
\right| \, ds \\
&&
{}+\Bigg| 4 \pi \!\int_0^t \ovR_m(s) \ovM_m(s) \Bigl[ 
\left( e^{(\l+\m)(s,\ovR_m(s))} - e^{(\ovla+\ovmu)(s,\ovR_m(s))}
\right) \j(s,\ovR_m(s)) \\
&&
\quad\quad\quad\quad\quad\quad\quad\quad
- \dot{\ovR}_m(s) \left( e^{2 \l(s,\ovR_m(s))}
- e^{2 \ovla(s,\ovR_m(s))} \right) \r(s,\ovR_m(s)) \Bigr] \, ds 
\Bigg|\\
&&
\quad+\Bigg| 4 \pi \!\int_0^t \ovR_m(s) \ovM_m(s) \Bigl[ 
e^{(\ovla+\ovmu)(s,\ovR_m(s))} \left( \j(s,\ovR_m(s))-\ovvj(s,\ovR_m(s))
\right) \\
&&
\quad\quad\quad\quad\quad\quad\quad\quad
-\dot{\ovR}_m(s) e^{2 \ovla(s,\ovR_m(s))} \left(
\r(s,\ovR_m(s))-\ovvrho(s,\ovR_m(s))\right) \Bigr] \, ds \Bigg|\\
&&
\quad+\Bigg| 4 \pi \int_0^t \ovR_m(s) \ovM_m(s) \Bigl[ 
e^{(\ovla+\ovmu)(s,\ovR_m(s))}\left(\ovvj(s,\ovR_m(s))-
\ovj(s,\ovR_m(s))\right) \\
&&
\quad\quad\quad\quad\quad\quad\quad\quad
-\dot{\ovR}_m(s) e^{2 \ovla(s,\ovR_m(s))} \left(
\ovvrho(s,\ovR_m(s))-\ovrho(s,\ovR_m(s))\right) \Bigr] \, ds \Bigg|.
\eeas
The troublesome term is the last one containing the differences 
$\ovvj-\ovj$ and
$\ovvrho-\ovrho$. If we denote it by $F_4$, then
\beas
|F_1| 
&\leq& 
C \e^3 \int_0^t \Bigl[ \|m(s)-\ovm(s)\|_\infty + \|\m(s)-\ovmu(s)\|_\infty
+ \|\l(s)-\ovla(s)\|_\infty\\
&&\quad\quad\quad\quad + \|\r(s)-\ovvrho(s)\|_\infty
        + \|\j(s)-\ovvj(s)\|_\infty\Bigr] \, ds + |F_4|.
\eeas
To continue we introduce some abbreviations:
\[
h_m(s):= \frac{e^{2 \ovla(s,\ovR_m(s))}}{\ovR_m(s)} \ovM_m(s),\
g_n(s) := e^{(\m-\l)(s,R_n(s))},\ 
\ovg_n(s) := e^{(\ovmu-\ovla)(s,R_n(s))}.
\]
The definitions of $g_{\ovn}$ and $\ovg_{\ovn}$ are analogous to those of
$g_n$ and $\ovg_n$ with $\ovR_n$ instead of $R_n$, and in particular
\[
\frac{W_n(s)}{E_n(s)} g_n(s) = \dot{R}_n(s),\ 
\frac{\ovW_n(s)}{\ovE_n(s)} \ovg_{\ovn}(s) = \dot{\ovR}_n(s).  
\]
Moreover,
$H^\d_{m,n}(s) := \chid(R_m(s)-R_n(s))$. The definitions of 
$H^\d_{\ovm,n}$, $H^\d_{\ovm,\ovn}$, and $H^\d_{m,\ovn}$ should be obvious.
Inserting these abbreviations and
the definition of the source terms into $F_4$ yields
\beas
|F_4| 
&=& 
\Bigg| \frac{1}{\d} \int_0^t \sum_{n=1}^N h_m(s) \Bigl[ \ovg_{\ovm}(s)
\Bigl( W_n M_n H^\d_{\ovm,n}
-\ovW_n \ovM_n H^\d_{\ovm,\ovn} \Bigr)(s) \\
&&
\quad\quad\quad\quad\quad\quad 
-\dot{\ovR}_m(s) 
\Bigl(E_n M_n H^\d_{\ovm,n}-\ovE_n \ovM_n H^\d_{\ovm,\ovn} \Bigl)(s) 
\Bigr]\,ds \Bigg|\\
&\leq&  
\Bigg| \frac{1}{\d} \int_0^t \sum_{n=1}^N h_m(s)
\ovW_n(s) \ovM_n(s) H^\d_{\ovm,n}(s) 
\bigl( \ovg_{\ovn} -\ovg_n \bigr)(s) \, ds \Bigg| \\
&&
{} +\Bigg| \frac{1}{\d} \int_0^t \sum_{n=1}^N h_m(s)
\ovg_{\ovm}(s) H^\d_{\ovm,n}(s) 
\bigl( W_n M_n - \ovW_n \ovM_n \bigr)(s) \, ds \Bigg| \\
&&
{} +\Bigg| \frac{1}{\d} \int_0^t \sum_{n=1}^N h_m(s)
\ovM_n(s)  H^\d_{\ovm,n}(s) 
\left( \ovg_n \ovW_n - g_n W_n \frac{\ovE_n}{E_n} \right)(s)\, ds \Bigg| \\
&&
{} +\Bigg| \frac{1}{\d} \int_0^t \sum_{n=1}^N h_m(s)
\dot{\ovR}_m(s) H^\d_{\ovm,n}(s) 
\bigl( \ovE_n \ovM_n - E_n M_n \bigr)(s) \, ds \Bigg| \\
&&
{}+ \Bigg| \frac{1}{\d} \int_0^t \sum_{n=1}^N h_m(s)
\ovW_n(s) \ovM_n(s) 
\bigl( \ovg_{\ovm} -\ovg_{\ovn}\bigr)(s) 
\bigl( H^\d_{\ovm,n} - H^\d_{\ovm,\ovn} \bigr)(s)
\, ds \Bigg| \\
&&
{} + \Bigg| \frac{1}{\d} \int_0^t \sum_{n=1}^N h_m(s)
\ovE_n(s) \ovM_n(s)\\
&&{}
\qquad\quad\quad\quad\quad
\Bigl( H^\d_{\ovm,n}(s) \bigl( \dot{R}_n -\dot{\ovR}_m \bigr)(s)
- H^\d_{\ovm,\ovn}(s) \bigl( \dot{\ovR}_n - \dot{\ovR}_m \bigr)(s) 
\Bigr)  ds \Bigg|\\
&=:& 
F_5 + F_6 + F_7 + F_8 +F_9 + F_{10}.
\eeas
Let
$I_1(s) := \{ n \in \CN \mid |R_n(s) - \ovR_m(s)| \leq \d\}$,
$I_2(s) := \{ n \in \CN \mid |\ovR_n(s) - \ovR_m(s)| \leq \d\}$. 
Then by Lemma~\ref{Icount}, $|I_1(s)| \leq C \e^{-3} \d$ and 
$\sum_{n \in I_2(s)} \ovM_n(s) \leq C \d$.
Hence,
\beas
|F_5| 
&\leq& 
C \frac{\e^6}{\d} \int_0^t \sum_{n \in I_1(s)} 
|R_n(s)\!-\!\ovR_n(s)| \, ds
\leq C \e^3 \int_0^t \NormRs \, ds,\\
|F_6| 
&\leq& 
C \frac{\e^6}{\d} \int_0^t \sum_{n \in I_1(s)} 
\Bigl( |W_n(s)-\!\ovW_n(s)|+\e^{-3}|M_n(s)-\!\ovM_n(s)|\Bigr) \, ds\\
&\leq& 
C \e^3 \int_0^t \Bigl( \NormWs + \NormMs \Bigr)\, ds,\\
|F_7| 
&\leq& 
C \frac{\e^6}{\d} \int_0^t \sum_{n \in I_1(s)} 
\Bigl( \|\m(s)-\ovmu(s)\|_\infty +\|\l(s)-\ovla(s)\|_\infty\\[-2mm]
&&
{}\quad\quad\quad\quad\quad\quad\quad+|R_n(s)-\ovR_n(s)|
+|W_n(s)-\ovW_n(s)| \Bigr)  ds\\
&\leq& 
C \e^3 \int_0^t \Bigl( \|\m(s)-\ovmu(s)\|_\infty +\|\l(s)-\ovla(s)\|_\infty\\
&&
{}\quad\quad\quad\quad+\NormRs+\NormWs \Bigr) \, ds,\\
|F_8| 
&\leq& 
C \frac{\e^6}{\d} \int_0^t \sum_{n \in I_1(s)}
\Bigl(|R_n(s)-\ovR_n(s)| + |W_n(s)-\ovW_n(s)|\\[-2mm]
&&
\qquad \qquad \qquad \qquad \qquad \qquad \qquad
+ \e^{-3}|M_n(s)-\ovM_n(s)|\Bigr) \, ds\\
&\leq& 
C \e^3 \int_0^t \Bigl( \NormRs + \NormWs + \NormMs \Bigr)\, ds,\\
|F_9| 
&\leq& 
C \frac{\e^3}{\d} \int_0^t \sum_{n \in I_1(s) \cup I_2(s)} 
\ovM_n(s)\,|\ovR_m(s)-\ovR_n(s)| \\[-2mm]
&&
\qquad \qquad \qquad \qquad \qquad \qquad
\min\left\{|R_n(s)-\ovR_n(s)|/\d,1\right\} \, ds\\
&\leq& 
C \frac{\e^3}{\d} \int_0^t 
\sum_{n \in I_1(s) \cup I_2(s)} \ovM_n(s)\, \left( \d + \NormRs \right)\\[-2mm]
&&
\qquad \qquad \qquad \qquad \qquad \qquad\quad
\min\left\{\NormRs/\d,1\right\} \, ds\\
&\leq& 
C \frac{\e^3}{\d} \int_0^t \NormRs 
\sum_{n \in I_1(s) \cup I_2(s)} \ovM_n(s) \, ds\\
&\leq& C \e^3 \int_0^t \NormRs \, ds.
\eeas
To deal with $F_{10}$ 
let $g_{n,m}(s) := h_m(s) \ovE_n(s) \ovM_n(s)$. 
By Lemma~\ref{discrbounds}, $|\dot{\ovM}_n| \leq C \e^3$, hence
for all $n,m \in \CN$,
\[ 
\|g_{n,m}\|_\infty \leq C \e^6,\quad
\|\dot g_{n,m}\|_\infty \leq C \e^6.
\]
Thus Lemma~\ref{gmn} yields
\beas
|F_{10}|
&\leq& 
C  \,\e^{-3}\max_{k,l\in\CN}
\bigl( \|g_{k,l}\|_\infty + \|\dot g_{k,l}\|_\infty\bigr)\,
\NormRt\\
&\leq& C \e^3 \NormRt.
\eeas
Putting all our estimates together and observing Lemma~\ref{mu-mub}, we get
\[
\frac{|\ovvM_m(t)-\ovM_m(t)|}{\e^3} \leq C \int_0^t 
\frac{|\ovvM_m(s)-\ovM_m(s)|}{\e^3} \, ds
+ C \int_0^t e(s) \, ds + C \NormRt. 
\]
Taking the maximum over $m \in \CN$ and using
a Gronwall argument completes the proof.
\prfe

\begin{lemma} \label{w-wb}
With $e(s)$ defined as in Lemma~\ref{mbb-mb},
\[ 
\NormWt \leq C \int_0^t e(s) \, ds + C \, \NormRt,\ t \in [0,\Ted].
\]
\end{lemma}
\prf
The proof is similar to the one of Lemma \ref{mbb-mb}.
Let $t \in [0,\Ted]$ and $m \in \CN$. From the differential equations
for $W_m$ and $\ovW_m$,
\beas 
W_m(t) - \ovW_m(t) 
&=&
\int_0^t \Biggl[\frac{ e^{(\m-\l)(s,R_m(s))} L_m }{R_m^3(s) E_m(s)}
- \dot\l(s,R_m(s))\, W_m(s)\\
&&
\qquad
- e^{(\m-\l)(s,R_m(s))} \m'(s,R_m(s))\, E_m(s)
- \frac{e^{(\ovmu-\ovla)(s,\ovR_m(s))} L_m}{\ovR_m^3(s) \ovE_m(s)} \\
&&
\quad
+ \ovlap(s,\ovR_m(s))\, \ovW_m(s)
+ e^{(\ovmu-\ovla)(s,\ovR_m(s))} 
\ovmus(s,\ovR_m(s))\, \ovE_m(s) \Biggr] ds.
\eeas
Inserting the equations for $\dot\l$, $\m'$, $\ovlap$ and $\ovmus$
yields
\beas
&&
W_m(t) - \ovW_m(t) = \\
&&
\int_0^t \biggl[\frac{e^{(\m-\l)(s,R_m(s))} L_m}{R_m^3(s) E_m(s)}
-\frac{e^{(\ovmu-\ovla)(s,\ovR_m(s))} L_m}{\ovR_m^3(s) \ovE_m(s)}\biggr]\, ds\\
&&
{}+\int_0^t \!\!\biggl[ e^{(\ovmu+\ovla)(s,\ovR_m(s))}
\ovE_m(s)\frac{\ovm(s,\ovR_m(s))}{\ovR_m^2(s)}  
- e^{(\m+\l)(s,R_m(s))}
E_m(s)\frac{m(s,R_m(s))}{R_m^2(s)} \biggr] ds\\
&&
{}+\int_0^t \biggl[ 4 \pi R_m(s) e^{(\m+\l)(s,R_m(s))} 
\Bigl( W_m(s) \j(s,R_m(s))- E_m(s) \,\, p(s,R_m(s)) \Bigr)\\
&&
\quad\quad\quad
-4 \pi \ovR_m(s) e^{(\ovmu+\ovla)(s,\ovR_m(s))} \Big( \ovW_m(s)
\ovj(s,\ovR_m(s)) - \ovE_m(s) \ovp(s,\ovR_m(s))\Big) \biggr] ds\\
&&
=: 
F_1 + F_2 + F_3;
\eeas
we have changed the order of terms to group the 
un-integrated source terms together.
Using the same calculations as in the lemmas above we find
\beas 
|F_1| 
&\leq& 
C \int_0^t \Bigl( \|\m(s)-\ovmu(s)\|_\infty +
\|\l(s)-\ovla(s)\|_\infty \\
&&
\quad\quad\quad\quad +\NormRs + \NormWs \Bigr)\, ds,\\
|F_2|
&\leq& 
C \int_0^t \Bigl( \|\m(s)-\ovmu(s)\|_\infty +
\|\l(s)-\ovla(s)\|_\infty + \|m(s) - \ovm(s)\|_\infty \\
&&
\qquad\qquad\qquad\qquad\qquad\quad +\NormRs + \NormWs \Bigr)\, ds,\\
|F_3| 
&\leq& 
4 \pi
\biggl|\int_0^t \biggl[R_m(s) e^{(\m+\l)(s,R_m(s))}
\Bigl( W_m(s) \j(s,R_m(s)) - E_m(s) p(s,R_m(s)) \Bigr)\\
&&
\quad\quad -\ovR_m(s) e^{(\ovmu+\ovla)(s,\ovR_m(s))} 
\Bigl( \ovW_m(s) \j(s,\ovR_m(s))- \ovE_m(s) p(s,\ovR_m(s))\Bigr) 
\biggr] ds \biggr| \\ 
&&
{}+ 4 \pi \biggl| \int_0^t \ovR_m(s) e^{(\ovmu+\ovla)(s,\ovR_m(s))}
\biggl[\ovW_m(s) \left(\j(s,\ovR_m(s)) - \ovvj(s,\ovR_m(s)) \right)\\
&&
\qquad\qquad\quad\quad\quad\quad\quad\quad\quad\quad\quad\quad
- \ovE_m(s) 
\left( p(s,\ovR_m(s))-\ovvp(s,\ovR_m(s)) \right)\biggr]\, ds \biggr| \\
&&
{}+ 4 \pi \biggl| \int_0^t  \ovR_m(s) 
e^{(\ovmu+\ovla)(s,\ovR_m(s))}
\biggl[\ovW_m(s) \left( \ovvj(s,\ovR_m(s)) - \ovj(s,\ovR_m(s)) \right)\\
&&
\qquad\qquad\quad\quad\quad\quad\quad\quad\quad\quad\quad\quad
- \ovE_m(s) 
\left( \ovvp(s,\ovR_m(s))-\ovp(s,\ovR_m(s)) \right)\biggr]\, ds \biggr|.
\eeas
Denote the last term---the one with the differences $\ovvj-\ovj$ and
$\ovvp-\ovp$---by $F_4$. Then 
\beas
|F_3| 
&\leq& 
C \int_0^t \biggl[ \|\m(s)-\ovmu(s)\|_\infty
+ \|\l(s)-\ovla(s)\|_\infty +\NormRs\\
&&
\quad\quad\quad +\NormWs+ \|p(s)-\ovvp(s)\|_\infty
+ \|\j(s)-\ovvj(s)\|_\infty\biggr] \, ds + |F_4|.
\eeas
We insert the definitions of the source terms into $F_4$, 
define $k_m(s):=e^{2 \ovla(s,\ovR_m(s))} \ovE_m(s)/\ovR_m(s)$, 
recall the definitions of $g_m$, 
$H^\d_{m,n}$ etc.\ from the proof of Lemma \ref{mbb-mb}
and note that
$\ovW_m \ovg_{\ovm} /\ovE_m  = \dot{\ovR}_m$.
Then
\beas
|F_4|
&=& 
\Biggl| \frac{1}{\d} \int_0^t \sum_{n=1}^N k_m(s) \ovg_{\ovm}(s)
\Biggl[ \frac{\ovW_m(s)}{\ovE_m(s)} 
\Bigl( W_n M_n H^\d_{\ovm,n} 
-\ovW_n \ovM_n H^\d_{\ovm,\ovn} \Bigr)(s) \\
&&
\qquad\qquad\quad\quad\quad\quad\quad\quad\quad\quad\quad
+\biggl(\frac{\ovW_n^2}{\ovE_n} \ovM_n  H^\d_{\ovm,\ovn} 
- \frac{W_n^2}{E_n} M_n  H^\d_{\ovm,n} \biggr)(s) \Biggr] ds \Biggr|\\
&\leq&  
\Biggl| \frac{1}{\d} \int_0^t \sum_{n=1}^N k_m(s)
\ovM_n(s) H^\d_{\ovm,n}(s) 
\left( \frac{\ovW_n W_n}{E_n} g_n - \frac{\ovW_n
\ovW_n}{\ovE_n} \ovg_{\ovn} \right)(s) ds \Biggr| \\
&&
\quad+\Biggl| \frac{1}{\d} \int_0^t \sum_{n=1}^N k_m(s)
\frac{\ovW_m(s)}{\ovE_m(s)} 
\ovg_{\ovm}(s) H^\d_{\ovm,n}(s) 
\bigl( W_n M_n - \ovW_n \ovM_n \bigr)(s) \, ds \Biggr| \\
&&
\quad+\Biggl| \frac{1}{\d} \int_0^t \sum_{n=1}^N k_m(s)
\ovg_{\ovm}(s) H^\d_{\ovm,n}(s) 
\left( \frac{\ovW_n^2}{\ovE_n} 
\ovM_n - \frac{W_n^2}{E_n} M_n \right)(s) \, ds \Biggr| \\
&&
\quad+\Bigg| \frac{1}{\d} \int_0^t \sum_{n=1}^N k_m(s)
\frac{\ovW_n^2(s)}{\ovE_n(s)} \ovM_n(s) 
\bigl( \ovg_{\ovm} - \ovg_{\ovn} \bigr)(s) 
\bigl( H^\d_{\ovm,\ovn} - H^\d_{\ovm,n} \bigr)(s) \, ds \Biggr| \\
&&
\quad+\Biggl| \frac{1}{\d} \int_0^t \sum_{n=1}^N k_m(s)
\ovW_n(s) \ovM_n(s)\\[-2mm]
&&
\qquad \qquad \qquad \quad \Bigl( H^\d_{\ovm,n}(s) 
\left(\dot{\ovR}_m - \dot{R}_n \right)(s) 
-H^\d_{\ovm,\ovn}(s) \left( \dot{\ovR}_m - \dot{\ovR}_n 
\right)(s) \Bigr)  ds \Biggr|\\
&=:& 
F_5 + F_6 + F_7 + F_8 + F_9.
\eeas
Let
$I_1(s) := \{ n \in \CN \mid |R_n(s)-\ovR_m(s)| \leq \d\},\
I_2(s) := \{ n \in \CN \mid |\ovR_n(s)-\ovR_m(s)| \leq \d\}$. 
Then by Lemma~\ref{Icount},
$|I_1(s)| \leq C \e^{-3} \d$ and 
$\sum_{n \in I_2(s)} \ovM_n(s) \leq C \d$. Hence,
\beas
|F_5| 
&\leq& 
C \int_0^t \Bigl( \|\m(s)-\ovmu(s)\|_\infty +\|\l(s)-\ovla(s)\|_\infty  \\
&&
\quad\quad\quad
+ \NormRs + \NormWs \Bigr) \, ds,\\
|F_6| 
&\leq& 
C \frac{\e^3}{\d} \int_0^t \sum_{n \in I_1(s)} 
\Bigl( |W_n(s)-\ovW_n(s)|+\e^{-3}|M_n(s)-\ovM_n(s)| \Bigr) \, ds\\
&\leq& 
C \int_0^t \Bigl( \NormWs + \NormMs \Bigr) \, ds,\\
|F_7| 
&\leq& 
C \frac{\e^3}{\d} \int_0^t \sum_{n \in I_1(s)}
\Bigl( |R_n-\ovR_n| + |W_n-\ovW_n| + \e^{-3}|M_n-\ovM_n|
\Bigr)(s) \, ds\\
&\leq& 
C \int_0^t \Bigl( \NormRs + \NormWs + \NormMs \Bigr)\, ds,\\
|F_8| 
&\leq& 
\frac{C}{\d} \int_0^t \sum_{n \in I_1(s)\cup I_2(s)} 
\ovM_n(s)\Bigl(|\ovR_m-\ovR_n| 
\min\left\{|R_n-\ovR_n|/\d, 1 \right\} \Bigr)(s) ds\\
&\leq& 
\frac{C}{\d} \int_0^t \sum_{n \in I_1(s)\cup I_2(s)} \ovM_n(s)
\Bigl( (\d+\|R-\ovR\|) 
\min\left\{\|R-\ovR\|/\d, 1 \right\} \Bigr)(s) ds\\
&\leq& 
C \int_0^t \NormRs \, ds.
\eeas
Finally, let $g_{n,m}(s) := k_m(s) \ovW_n(s) \ovM_n(s)$. 
By Lemma~\ref{discrbounds},  
$|\dot{\ovM}_n(s)| \leq C \e^3$, hence
\[ 
\|g_{n,m}\|_\infty\leq C \e^3,\
\|\dot g_{n,m}\|_\infty\leq C \e^3,\ n,m \in \CN.
\]
Thus Lemma~\ref{gmn} yields
\[
|F_9| \leq C \,\e^{-3} \max_{k,l \in \CN}
\bigl( \|g_{k,l}\|_\infty+\|\dot g_{k,l}\|_\infty\bigr)\NormRt
\leq C \NormRt.
\]
Collecting all these estimates and recalling Lemma~\ref{mu-mub}
completes the proof.
\prfe
 
We are finally ready to prove our main results:\\
\smallskip\noindent
{\bf Proof of Theorem~\ref{main} and Theorem~\ref{mainw}.}\\
By Lemma~\ref{mbb-mb}, 
\beas
\|\ovM(t)-\ovvM(t)\| 
&\leq& 
C \Bigl( \NormRt + \Delta\Bigr) \\
&&
{} + C \int_0^t \Bigl(\NormWs +\NormMs \Bigr) ds.
\eeas
where 
\[
\Delta := \sup_{t \in [0,\Ted]} 
\Bigl( \|\r(t)-\ovvrho(t)\|_\infty + \|p(t)-\ovvp(t)\|_\infty +
\|\j(t)-\ovvj(t)\|_\infty\Bigr).
\]
Combining this with Lemma \ref{m-mbb} yields
\beas
\NormMt 
&\leq& 
C 
\Bigl( \NormRt + \Delta + \e\Bigr)\\
&&
{} + C \int_0^t \Bigl(\NormWs +\NormMs \Bigr) ds.
\eeas
Next, by Lemma~\ref{w-wb},
\beas
\NormWt 
&\leq& 
C \Bigl( \NormRt + \Delta\Bigr) \\
&&
{} + \int_0^t \Bigl(\NormWs +\NormMs \Bigr) ds.
\eeas
Adding the last two estimates and applying
a Gronwall argument implies
\be \label{wmrest}
\NormWt+\NormMt\leq C \Bigl( \NormRt + \Delta + \e \Bigr).
\ee
Using Lemma~\ref{dr-drb} and Lemma~\ref{mu-mub} we find
\beas
|R_n(t)-\ovR_n(t)|
&\leq& 
C \int_0^t \Bigl( |R_n(s)-\ovR_n(s)|+|W_n(s)-\ovW_n(s)|\\
&&
\qquad\qquad+\|\m(s)-\ovmu(s)\|_\infty +\|\l(s)-\ovla(s)\|_\infty \Bigr) \, ds\\
&\leq& 
C \Delta \\
&&
{}+ C \int_0^t \Bigl(\NormRs  + \NormWs + \NormMs \Bigr)\, ds
\eeas
for all $t \in [0,\Ted]$ and $n \in \CN$. Taking the maximum over 
$n \in \CN$ and inserting (\ref{wmrest}) yields
\[ 
\NormRt \leq C \left(\Delta + \e + \int_0^t \NormRs\, ds\right)
\]
which via Gronwall implies that
\be \label{rdiffest}
\NormRt \leq C \left( \Delta + \e \right).
\ee
Now we apply Lemma~\ref{rho-rhobb}:
\be \label{Dest}
\Delta \leq C
\left\{
\begin{array}{cl}
\d + \e/\d &\ \mbox{in case of Theorem~\ref{mainw}},\\
\d^2 + \e/\d &\ \mbox{in case of Theorem~\ref{main}}.
\end{array} \right.
\ee
Inserting this into (\ref{rdiffest})
proves the corresponding error estimates on $\NormRt$ in the two theorems. 
The error estimates for $\NormWt$ and $\NormMt$ follow from (\ref{wmrest}),
and those for $\|\m(t)-\ovmu(t)\|_\infty$, 
$\|\l(t)-\ovla(t)\|_\infty$, and $\|m(t)-\ovm(t)\|_\infty$ follow
from Lemma~\ref{mu-mub}.
Under the assumption of Theorem~\ref{main}
the estimates for $\|\r(t)-\ovrho(t)\|_\infty$,
$\|p(t)-\ovp(t)\|_\infty$ and $\|\j(t)-\ovj(t)\|_\infty$
are valid due to Lemma~\ref{rho-rhobb} and
Lemma~\ref{rhob-rhobb}.

It remains to establish the assertions on the length $\Ted$ of our
approximation interval. Let $C$ denote
a constant for which the error estimates in Theorem~\ref{main}
or Theorem~\ref{mainw} hold;
note that $C$ is independent of $\d$ and $\e$. For Theorem~\ref{mainw},
let $\d$ be so small that
\[ 
\d + \frac{\e}{\d} \leq \frac{1}{32 D^2 C}. 
\]
Then
\beas
\ovR_n(t) 
&\leq& 
R_n(t) + \NormRt 
\leq D + C \left( \d + \frac{\e}{\d} \right) 
\leq \frac{3}{2} D,\\
\frac{1}{\ovR_n(t)} 
&=& 
\frac{1}{R_n(t)} + \frac{R_n(t)-\ovR_n(t)}
{R_n(t) \ovR_n(t)} 
\leq D + 2 D^2 C \left( \d + \frac{\e}{\d} \right) 
\leq \frac{3}{2} D,\\
|\ovW_n(t)| 
&\leq& 
|W_n(t)| + \NormWt
\leq D + C \left( \d + \frac{\e}{\d} \right) 
\leq \frac{3}{2} D,\\
|\ovM_n(t)| 
&\leq& 
|M_n(t)| + \e^3 \NormMt
\leq \e^3 \left( D + C \left( \d + \frac{\e}{\d} \right) \right) 
\leq \e^3 \frac{3}{2} D,\\
e^{2\ovla(t,r)}
&=& 
e^{2\l(t,r)} + \frac{2}{r} (\ovm - m)(t,r) e^{2(\l+\ovla)(t,r)}
\leq  D  + 16 D^3 C \left(\d + \frac{\e}{\d} \right) 
\leq  \frac{3}{2} D
\eeas
for all $t \in [0,\Ted]$ and $r > 0$ since $\ovm(t,r)=m(t,r)=0$ 
for $r < 1/(4D)$. 
Thus $\Ted=T$ since otherwise the approximation
interval could be extended beyond $\Ted$. 

For Theorem~\ref{main}, let $\d$ be so small that
\[ 
\d + \frac{\e}{\d} \leq \frac{1}{32 D^2 C} \ \wedge \ 
\d + \frac{\e}{\d^2} \leq \frac{D}{2 C}.
\]
Then the above inequalities hold for all $t \in [0,\Ted]$, and 
in addition
\[ 
\|\ovrho(t)\|_\infty \leq \|\r(t)\|_\infty+
\|\ovrho(t)-\r(t)\|_\infty\leq D + C 
\left( \d + \frac{\e}{\d^2} \right) \leq \frac{3}{2} D.
\]
With the same argument as above, $\Ted=T$.
This completes the proof of the main results.
     
\section{The fully discretized approximation}
\setcounter{equation}{0}

In setting up a discretized version of a system like the Vlasov-Einstein
system discretizing the phase space is the major step.
For example, in several of the papers on numerical schemes for
the Vlasov-Poisson or Vlasov-Maxwell systems mentioned in the
introduction only this step is analyzed.
However, due to the 
presence of the un-integrated source terms in the characteristic 
system, we can not discretize our system in time as easily . 
We have to modify the evolution equations slightly in order to 
get a convergence result, the idea being to discretize after 
integrating by parts in Lemma~\ref{gmn}. Recall that the latter maneuver
was essential for the convergence proof of the semi-discretized
approximation, and we have to set up things in such a way that there
is a substitute for this in the fully discretized case. 

Throughout this section the assumptions of
Theorem \ref{main} are to be satisfied. Moreover, we assume
$\e \leq \d^2$ and $\d \leq 1/(8D)$.
Let $0 < \tau \ll 1$ be a small parameter, the time step. 
Before we discretize the evolution equations for 
$\ovR_n(t),\ \ovW_n(t),\ \ovM_n(t)$
we regroup the terms in these equations as follows;
only the equations for $\ovW_n(t)$ and $\ovM_n(t)$ need to be considered,
the discretization of the equation for $\ovR_n(t)$ being trivial:
\beas
\dot{\ovW}_n(t)
&=& 
e^{(\ovmu-\ovla)(t,\ovR_n(t))} \frac{L_n}{\ovR_n^3(t) \ovE_n(t)}
- e^{(\ovmu+\ovla)(t,\ovR_n(t))} \frac{\ovm(t,\ovR_n(t))}{\ovR_n^2(t)} 
\ovE_n(t)\\
&+&
\frac{e^{2 \ovla(t,\ovR_n(t))}}{\ovR_n(t)}
4 \pi \ovR_n^2(t) e^{(\ovmu-\ovla)(t,\ovR_n(t))}
\left(\ovW_n(t)\ovj(t,\ovR_n(t))  - \ovE_n(t) \ovp(t,\ovR_n(t)) \right)\\
&=& 
e^{(\ovmu-\ovla)(t,\ovR_n(t))} \frac{L_n}{\ovR_n^3(t) \ovE_n(t)}
- e^{(\ovmu+\ovla)(t,\ovR_n(t))} \frac{\ovm(t,\ovR_n(t))}{\ovR_n^2(t)} 
\ovE_n(t)\\
&&
{}+ \frac{e^{2 \ovla(t,\ovR_n(t))}}{\ovR_n(t)}
\frac{1}{\d} \sum_{m=1}^N \ovM_m(t) \ovW_m(t) 
\ovE_n(t) 
\frac{d}{dt}\chi(\ovR_n(t)-\ovR_m(t))\\
&&
{}+\frac{e^{2 \ovla(t,\ovR_n(t))}}{\ovR_n(t)}
\frac{1}{\d} \sum_{m=1}^N \ovM_m(t) \ovW_m^2(t) 
\frac{\ovE_n(t)}{\ovE_m(t)} \chid(\ovR_n(t)-\ovR_m(t))\\ 
&&
\qquad \qquad \qquad \qquad \qquad
\left(e^{(\ovmu-\ovla)(t,\ovR_m(t))} - e^{(\ovmu-\ovla)(t,\ovR_n(t))}\right)\\
&=:& 
\ovfwan(t) + \ovfwbn(t) + \ovfwcn(t);
\eeas
note that
\beas
\frac{d}{dt}\chi(\ovR_n(t)-\ovR_m(t))
&=& 
\frac{1}{\d}\chid(\ovR_n(t)-\ovR_m(t))\\
&&
{} \qquad
\left(e^{(\ovmu-\ovla)(t,\ovR_n(t))}\frac{\ovW_n(t)}{\ovE_n(t)}
-e^{(\ovmu-\ovla)(t,\ovR_m(t))}\frac{\ovW_m(t)}{\ovE_m(t)}\right).
\eeas
Similarly,
\beas
\dot{\ovM}_n(t) 
&=& 
\ovM_n(t) \Biggl[ e^{(\ovmu+\ovla)(t,\ovR_n(t))} 
\frac{\ovm(t,\ovR_n(t))}{\ovR_n^2(t)} \frac{\ovW_n(t)}{\ovE_n(t)}\\
&&
{}+ \frac{e^{2 \ovla(t,\ovR_n(t))}}{\ovR_n(t)} 
4 \pi \ovR_n^2(t) e^{(\ovmu-\ovla)(t,\ovR_n(t))}
\left( \ovj(t,\ovR_n(t)) - 
\frac{\ovW_n(t)}{\ovE_n(t)}\ovrho(t,\ovR_n(t)) \right) \Biggr]\\
&=&
\ovM_n(t) \Biggl[e^{(\ovmu+\ovla)(t,\ovR_n(t))} 
\frac{\ovm(t,\ovR_n(t))}{\ovR_n^2(t)} \frac{\ovW_n(t)}{\ovE_n(t)}\\
&&
\qquad \qquad - \frac{e^{2 \ovla(t,\ovR_n(t))}}{\ovR_n(t)}
\frac{1}{\d} \sum_{m=1}^N \ovM_m(t) 
\ovE_m(t) \frac{d}{dt}\chi(\ovR_n(t)-\ovR_m(t)) \\
&&
\qquad \qquad +\frac{e^{2 \ovla(t,\ovR_n(t))}}{\ovR_n(t)} 
\frac{1}{\d} \sum_{m=1}^N \ovM_m(t) \ovW_m(t)\,\chid(\ovR_n(t)-\ovR_m(t))\\
[-1mm]
&&
\qquad \qquad \qquad \qquad\qquad \qquad\qquad
\left(e^{(\ovmu+\ovla)(t,\ovR_n(t))} -e^{(\ovmu+\ovla)(t,\ovR_m(t))} \right) 
 \Biggr]\\
&=:& 
\ovM_n(t) \Biggl[ \ovfman(t) + \ovfmbn(t) + \ovfmcn(t)\Biggr].
\eeas
The terms $\ovfwan(t),\ \ovfwbn(t),\ \ovfwcn(t)$ and  
$\ovfman(t),\ \ovfmbn(t),\ \ovfmcn(t)$ 
have natural discretizations given below;
note that $\frac{d}{dt}\chi(\ovR_n(t)-\ovR_m(t))$
is much better to discretize and analyze than the term 
$\frac{1}{\d} \left( \dot{\ovR}_n-\dot{\ovR}_m \right)
\chid(\ovR_n(t)-\ovR_m(t))$.
Discretizing the equations above leads to
the following Euler-like scheme:
\beas
\ovR_{n,i+1} 
&:=& 
\ovRni + \tau e^{(\ovmu_i-\ovla_i)(\ovRni)} \frac{\ovWni}{\ovE_{n,i}},\\
\ovW_{n,i+1} 
&:=& 
\ovWni + \tau \left( \ovfwani + \frac{e^{2 \ovla_i(\ovRni)}}{\ovRni} \ovfwbni
+ \ovfwcni  \right), \\
\ovM_{n,i+1} 
&:=& 
\ovMni + \tau \ovMni \left( \ovfmani  + \frac{e^{2 \ovla_i(\ovRni)}}{\ovRni} 
\ovfmbni + \ovfmcni \right),\\
\ovfwani 
&:=& 
e^{(\ovmu_i-\ovla_i)(\ovRni)} \frac{L_n}{\ovRni^3 \ovE_{n,i}}
- e^{(\ovmu_i+\ovla_i)(\ovRni)} \frac{\ovm_i(\ovRni)}{\ovRni^2} \ovE_{n,i},\\
\ovfwbni 
&:=& 
\frac{e^{2 \ovla_i(\ovRni)}}{\ovRni}
\sum_{m=1}^N \ovMmi \ovWmi \ovE_{n,i}\frac{\chi(\ovR_{n,i+1}-\ovR_{m,i+1})
- \chi(\ovR_{n,i}-\ovR_{m,i})}{\tau} ,\\
\ovfwcni 
&:=& 
\frac{e^{2 \ovla_i(\ovRni)}}{\ovRni} \frac{1}{\d}
\sum_{m=1}^N \ovMmi \ovWmi^2 \frac{\ovE_{n,i}}{\ovE_{m,i}}
\chid(\ovR_{n,i}-\ovR_{m,i})\\[-1mm]
&&
\qquad\qquad\qquad\qquad
\left(e^{(\ovmu_i-\ovla_i)(\ovRmi)}-e^{(\ovmu_i-\ovla_i)(\ovRni)} \right),\\
\ovfmani 
&:=& 
e^{(\ovmu_i+\ovla_i)(\ovRni)} \frac{\ovm_i(\ovRni)}{\ovRni^2} 
\frac{\ovWni}{\ovE_{n,i}},\\
\ovfmbni 
&:=& 
\frac{e^{2 \ovla_i(\ovRni)}}{\ovRni}
\sum_{m=1}^N \ovMmi \ovE_{m,i} \frac{\chi(\ovR_{n,i+1}-\ovR_{m,i+1})
- \chi(\ovR_{n,i}-\ovR_{m,i})}{\tau},\\
\ovfmcni 
&:=& 
\frac{e^{2 \ovla_i(\ovRni)}}{\ovRni} \frac{1}{\d}
\sum_{m=1}^N \ovMmi \ovWmi\, \chid(\ovR_{n,i}-\ovR_{m,i})\\ [-1mm]
&&
\qquad\qquad\qquad\qquad 
\left(e^{(\ovmu_i-\ovla_i)(\ovRni)}-e^{(\ovmu_i-\ovla_i)(\ovRmi)} \right),
\eeas
with initial data 
\[
(\ovR_{n,0},\ovW_{n,0},\ovM_{n,0}) := (\ovR_n(0),\ovW_n(0),\ovM_n(0)),
\ n \in \CN.
\]
We have abbreviated 
\[
\ovE_{n,i}:=\sqrt{1+\ovW^2_{n,i} + L_n/\ovR^2_{n,i}}\, ; 
\]
note that $\ovR_{n,i+1}$ can be computed without knowing
the terms $\ovfwani$ and $\ovfmani$ where $\ovR_{n,i+1}$ appears
in the difference quotient.
The definitions of the approximating source terms etc.\ 
at the $i$-th time step are completely analogous to the definitions of 
these quantities at time $t=i\tau$ in case of the phase space discretization.
For example, for $r > 0$, 
\[
\ovrho_i(r) := \frac{1}{4 \pi r^2 \d} \sum_{n=1}^N \ovMni \ovE_{n,i}
\chid(r-\ovRni),
\]
and analogously for $\ovp_i$ and $\ovj_i$. Inserting these
into the formulas for $\ovm(t)$, $\ovla(t)$, and $\ovmu(t)$
defines corresponding approximations at the $i$-th time step
of these quantities. 

To measure the errors for the fully discretized system
we define
\beas
e^R_i 
&:=& 
\max_{n \in \CN, k \in \{0,1,\dots,i\}} 
\left| \ovR_n(t_k)-\ovRnk \right|,\\
e^W_i 
&:=& 
\max_{n \in \CN, k \in \{0,1,\dots,i\}} 
\left| \ovW_n(t_k)-\ovWnk \right|,\\
e^M_i 
&:=& 
\max_{n \in \CN, k \in \{0,1,\dots,i\}} 
\e^{-3} \left| \ovM_n(t_k)-\ovMnk \right|;
\eeas
note that these are monotonous in $i$. The convergence result
for the fully discretized scheme is as follows:

\begin{theorem} \label{mainfudi}
Suppose the assumptions of Theorem~\ref{main} hold,
$\e \leq \d^2$, and let $\d$ and $\tau$ be sufficiently small. 
Then 
\beas
e^R_i + e^W_i + e^M_i 
&\leq& C \tau,\\
\|\ovm(i\tau)-\ovm_i\|_\infty + \|\ovmu(i\tau)-\ovmu_i\|_\infty
+ \|\ovla(i\tau)-\ovla_i\|_\infty 
&\leq& 
C \tau,\\
\|\ovrho(i\tau)-\ovrho_i\|_\infty + \|\ovp(i\tau)-\ovp_i\|_\infty
+ \|\ovj(i\tau)-\ovj_i\|_\infty 
&\leq& 
C \frac{\tau}{\d}
\eeas
for all $i\in \N$ such that $i\tau \in[0,T]$, the time interval
on which the solution $f$ is to be approximated.
\end{theorem}

As to the proof only some indications are given
since once the fully discretized system is set
up in the appropriate way, what remains is to establish a series of 
auxiliary results most of which have analogues in the semi-discretized 
approximation which were proven in that situation. The detailed proof
can be found in \cite{Ro}.
Let
\[
\Jed := 
\Bigl\{ i \in \N \mid i\tau \in [0,\Ted] \Bigr\}.
\]
The approximations are calculated as long as $i\in\Jed$ and
for all $n \in \CN$ and $k \in\{0,\dots,i\}$ the estimates
\[
\e^{-3} \ovMnk,\; |\ovWnk|,\; \ovRnk,\; \frac{1}{\ovRnk},\;
\|e^{2 \ovla_k}\|_\infty,\; \|\ovrho_k\|_\infty
\leq 4 D
\]
hold; observe the analogy to (\ref{teda}), (\ref{tedb}).
The set of all these indices $i$ is denoted by $\Jedt$
For these $i$ the error estimates from Theorem~\ref{mainfudi} are
established. This is then used to show that if $i\in\Jedt$ and $i+1 \in \Jed$
then $i+1 \in \Jedt$, provided the time step $\tau$ is sufficiently small.
Since $\Ted=T$ for sufficiently small $\d$ this completes this
indication of the proof.

\bigskip

{\bf Acknowledgment:} GR acknowledges support by the Wittgenstein 2000
Award of P.~A.~Markowich, funded by the Austrian FWF.


\begin{thebibliography}{10}

\bibitem{A1}
{\sc H.~Andr\'{e}asson}.
The Einstein-Vlasov system/Kinetic theory.
{\em Living Reviews in Relativity},
to appear

\bibitem{B}
{\sc C.~K.~Birdsall, A.~B.~Langdon}.
{\em Plasma Physics via Computer Simulation}.
McGraw-Hill (1985).

\bibitem{D}
{\sc J.~M.~Dawson}. 
Particle simulation of plasmas.
{\em  Rev. Modern Phys.} 
{\bf 55} (1983), 403--447.

\bibitem{V4}
{\sc K.~Ganguly, J.~T.~Lee, H.~D.~Victory,~Jr.}
On simulation methods for
Vlasov-Poisson systems with particles initially asymptotically distributed.
{\em SIAM J. Numer. Anal.} 
{\bf 28} (1991), 1574--1609.

\bibitem{V1}
{\sc K.~Ganguly, H.~D.~Victory,~Jr.} 
On the convergence of particle methods for
multidimensional Vlasov-Poisson systems.
{\em  SIAM J. Numer. Anal.} 
{\bf 26} (1989), 249--288 .

\bibitem{GS}
{\sc R.~Glassey, J.~Schaeffer}. 
Convergence of a particle method for the
relativistic Vlasov-Maxwell system.
{\em  SIAM J. Numer. Anal.} 
{\bf 28}, 1-25 (1991).

\bibitem{LP}
{\sc P.-L.~Lions, B.~Perthame}.
Propagation of moments and regularity for the 3-dimensional 
Vlasov-Poisson system.
{\em Invent.\ Math}.\ 
{\bf 105} (1991), 415--430.

\bibitem{OC}
{\sc I.~Olabarrieta, M.~W.~Choptuik}.
Critical phenomena at the threshold of black hole 
formation for collisionless matter in spherical symmetry.
{\em Phys.\ Rev. D.}\
{\bf 65} (2002), 024007

\bibitem{Pf}
{\sc K.~Pfaffelmoser}. 
Global classical solutions of the Vlasov-Poisson system in three dimensions 
for general initial data.
{\em J.\ Diff.\ Eqns}.\ 
{\bf 95} (1992), 281--303.

\bibitem{R1}
{\sc G.~Rein}.
Static solutions of the spherically symmetric Vlasov-Einstein system.
{\em Math.\ Proc.\ Camb.\ Phil.\ Soc.}\ 
{\bf 115} (1994), 559--570.


\bibitem{RH}
{\sc G. Rein}. 
{\em The Vlasov-Einstein System with Surface Symmetry.}
Habilitationsschrift, Munich, 1995.


\bibitem{R2}
{\sc G.~Rein}.
Static shells for the Vlasov-Poisson and Vlasov-Einstein systems.
{\em Indiana University Math.\ J.}\ 
{\bf 48} (1999), 335--346. 

\bibitem{RR1}
{\sc G.~Rein, A.~D.~Rendall}.
Global existence of solutions of the spherically symmetric
Vlasov-Einstein system with small initial data.
{\em Commun.\ Math.\ Phys.} 
{\bf 150} (1992), 561--583.

\bibitem{RR2}
{\sc G.~Rein, A.~D.~Rendall}.
The Newtonian limit of the spherically symmetric Vlasov-Einstein system.
{\em Commun.\ Math.\ Phys}.\ 
{\bf 150} (1992), 585--591 

\bibitem{RR3}
{\sc G.~Rein, A.~D.~Rendall}.
Smooth static solutions of the spherically symmetric Vlasov-Einstein system.
{\em Ann.\ de l'Inst.\ H.\ Poincar\'e, Physique Th\'eorique} 
{\bf 59} (1993), 383--397.

\bibitem{RR4}
{\sc G.~Rein, A.~D.~Rendall}.
Compact support of spherically symmetric equilibria in 
non-relativistic and relativistic galactic dynamics.
{\em Math.\ Proc.\ Camb.\ Phil.\ Soc.}\ 
{\bf 128} (2000), 363--380.

\bibitem{RRS1}
{\sc G.~Rein, A.~D.~Rendall, J.~Schaeffer}.
A regularity theorem for the spherically symmetric Vlasov-Einstein
system. 
{\em Commun.\ Math.\ Phys.} 
{\bf 168} (1995), 467--478.

\bibitem{RRS2}
{\sc G. Rein, A. Rendall, J. Schaeffer}. 
Critical collapse of collisionless matter:
A numerical investigation.
{\em  Phys. Rev. D} {\bf  58} (1998), 044007.

\bibitem{Ren1}
{\sc A.~D.~Rendall}.
Cosmic censorship and the Vlasov equation. 
{\em Class.\ Quantum Grav}. {\bf 9} (1992), L99--L104 .

\bibitem{Ren2}
{\sc A.~D.~Rendall}.
The Newtonian limit for asymptotically flat solutions 
of the Vlasov-Einstein system.
{\em Commun.\ Math.\ Phys}.\ 
{\bf  163} (1994), 89--112 

\bibitem{Ro}
{\sc T.~Rodewis}.
{\em Numerical treatment of the symmetric Vlasov-Poisson 
and Vlasov-Einstein system by particle methods}.
PhD thesis, Munich 1999.


\bibitem{Sch1}
{\sc J. Schaeffer}. 
Discrete Approximation of the Poisson-Vlasov System.
{\em Q.~Appl. Math.} 
{\bf 45} (1987), 59--73.

\bibitem{Sch2} 
{\sc J.~Schaeffer}. 
Global existence of smooth solutions to the Vlasov-Poisson system 
in three dimensions.
{\em Commun.\ Part.\ Diff.\ Eqns.}\ 
{\bf 16} (1991), 1313--1335.

\bibitem{Sh1}
{\sc S.~L.~Shapiro, S.~A.~Teukolsky}. 
Relativistic Stellar Dynamics on the computer.
I Motivation and numerical method.
{\em  Astrophys. J.} 
{\bf 298} (1985), 34--57.

\bibitem{Sh2}
{\sc S.~L.~Shapiro, S.~A.~Teukolsky}.
Relativistic Stellar Dynamics on the computer.
II Physical applications.
{\em  Astrophys. J.} 
{\bf 298} (1985), 58--79.

\bibitem{Sh3}
{\sc S.~L.~Shapiro, S.~A.~Teukolsky}.
Relativistic Stellar Dynamics on the computer. III.
{\em  Astrophys. J. (Letters)} 
{\bf 292} (1985), L~41.

\bibitem{Sh4}
{\sc S.~L.~Shapiro, S.~A.~Teukolsky}.
Relativistic Stellar Dynamics on the computer.
IV Collapse of a stellar cluster to a black hole.
{\em  Astrophys. J.} 
{\bf 307} (1986), 575--592.

\bibitem{V2}
{\sc H.~D.~Victory,~Jr., G.~Tucker, K.~Ganguly}.
The convergence analysis of fully discretized particle methods for solving
Vlasov-Poisson systems.
{\em  SIAM J. Numer. Anal.} 
{\bf 28} (1991), 955--989.

\bibitem{V3}
{\sc H.~D.~Victory,~Jr., E.~J.~Allen}.
The convergence theory of particle-in-cell 
methods for multidimensional Vlasov-Poisson systems.
{\em SIAM J. Numer. Anal.} 
{\bf 28} (1991), 1207--1241.


\end{thebibliography}
\end{document}